\documentclass{amsart}

\usepackage[mathcal]{eucal}
\usepackage{graphicx}
\usepackage{pinlabel}
\usepackage{amscd}

\newtheorem{theorem}{Theorem}
\newtheorem{lemma}[theorem]{Lemma}
\newtheorem{definition}[theorem]{Definition}
\newtheorem{prop}[theorem]{Proposition}
\newtheorem{corollary}[theorem]{Corollary}

\begin{document}

\title{Quantum Teichm\"uller space and Kashaev algebra}

\author{Ren Guo} 

\address{School of Mathematics, University of Minnesota, Minneapolis, MN, 55455}

\email{guoxx170@math.umn.edu}

\author{Xiaobo Liu}

\address{Department of Mathematics, Columbia University, 2990 Broadway, New York, NY 10027}

\email{xiaoboli@math.columbia.edu}

\thanks{}



\begin{abstract}
Kashaev algebra associated to a surface is a noncommutative deformation of the algebra of rational functions of Kashaev coordinates. For two arbitrary complex numbers, there is a generalized Kashaev algebra. The relationship between the shear coordinates and Kashaev coordinates induces a natural relationship between the quantum Teichm\"uller space and the generalized Kashaev algebra.
\end{abstract}

\maketitle

\section{Introduction}

A quantization of the Teichm\"uller space $\mathcal T(S)$ of a punctured surface $S$ was developed by Chekhov and Fock \cite{CF,Fo,FC} and, independently, by Kashaev \cite{Kas1, Kas2, Kas3, Kas4}. This is a deformation of the $\mathrm C^*$--algebra of functions on Teichm\"uller space $\mathcal T(S)$. The quantization was expressed in terms of self-adjoint operators on Hilbert spaces and the quantum dilogarithm function. Although these two approaches of quantization use the same ingredients, the relationship between them is still mysterious. Chekhov and Fock worked with shear coordinates of Teichm\"uller space while Kashaev worked with a new coordinate.

The pure algebraic foundation of Chekhov-Fock's quantization was established in \cite{Liu1} (see also \cite{BonLiu,BBL}). In this paper we investigate the algebraic aspect of Kashaev's quantization and establish a natural relationship between these two algebraic theories. This algebraic relationship should shed light on the two approach of operator-theoritical quantization of Teichm\"uller space.

\subsection{Quantum Teichm\"uller space}

Let's review the finite dimensional Chekhov-Fock's quantization following \cite{Liu1}.
Let $S$ be an oreinted surface of finite topological type, with genus $g$ and with $p\geq 1$ punctures, obtained by removing $p$ points $\{v_1,\ldots,v_p\}$
from a closed oriented surface $\bar S$ of genus $g$. If the Euler characteristic of $S$ is negative, i.e., $m:=2g-2+p>0,$ $S$ admits complete hyperbolic metrics. The \emph{Teichm\"uller space} $\mathcal T(S)$ of $S$ consists of all isotopy classes of complete hyperbolic metrics on $S$.

An \emph{ideal triangulation} of $S$ is a triangulation of the closed surface $\bar{S}$ whose vertex set is exactly
$\{v_1,\ldots,v_p\}$. Under a complete hyperbolic metric, an ideal triangulation of $S$ is realized as a proper 1--dimensional submanifold  whose complementary regions are hyperbolic ideal triangles. William Thurston \cite{Th} associated to each ideal triangulation a global coordinate system which is called \emph{shear coordinate} (see also \cite{Bon, Fo}). Given two ideal triangulations $\lambda$ and $\lambda'$, the corresponding coordinate changes are rational, so that there is a well-defined notion of rational functions on $\mathcal T(S)$.

For an ideal triangulation $\lambda$ and a number $q= \mathrm e ^{\pi \mathrm i \hbar} \in \mathbb C$, the \emph{Chekhov-Fock algebra} $\mathcal T_\lambda^q$ is the algebra over $\mathbb C$ defined by generators $X_1^{\pm1}$, $X_2^{\pm1}$, \dots, $X_{3m}^{\pm1}$ associated to the components of $\lambda$ and by relations $X_iX_j=q^{2\sigma^\lambda_{ij}}X_jX_i$, where the numbers $\sigma^\lambda_{ij}$ are integers determined by the combinatorics of the ideal triangulation $\lambda$. This algebra has a well-defined fraction division algebra $\widehat{\mathcal T}_\lambda^q$.

As one moves from one ideal triangulation $\lambda$  to another $\lambda'$, Chekhov and Fock \cite{Fo, FC, CF} (see also \cite{Liu1}) introduce \emph{coordinate change isomorphisms} $\Phi_{\lambda\lambda'}^q: \widehat{\mathcal T}_{\lambda'} ^q \rightarrow \widehat{\mathcal T}_\lambda^q$ which satisfy the natural property that $\Phi_{\lambda''\lambda'}^q \circ \Phi_{\lambda'\lambda}^q = \Phi_{\lambda''\lambda}^q $ for any ideal triangulations $\lambda$, $\lambda'$, $\lambda''$. In a triangulation independent way, this associates to the surface
$S$ the algebra $\widehat{\mathcal T}_S^q$ defined as the quotient of the family of all $\widehat{\mathcal T}_\lambda^q$, with $\lambda$ ranging over ideal triangulations of the surface $S$, by the equivalence relation
that identifies $\widehat{\mathcal T}_\lambda^q$ and $\widehat{\mathcal T}_{\lambda'}^q $ by the coordinate change
isomorphism $\Phi_{\lambda\lambda'}^q$. The algebra $\widehat{\mathcal T}_S^q$ is called the \emph{quantum Teichm\"uller space} of the surface $S$. It turns out that $\Phi_{\lambda\lambda'}^1$ is just
the corresponding shear coordinate changes. Therefore, the quantum Teichm\"uller space $\widehat{\mathcal T}_S^q$
is a noncommutative deformation of the algebra of rational functions on the
Teichm\"uller space $\mathcal T(S)$.

\subsection{Generalized Kashaev algebra}

A \emph{decorated ideal triangulation} of a punctured surface $S$ is an ideal triangulation such that the ideal triangles are numerated and there is a mark at a corner of each triangle. Kashaev \cite{Kas1} introduced a new coordinate associated to a decorated ideal triangulation of $S$. A Kashaev coordinate associated to a decorated ideal triangulation is a vector in $\mathbb{R}^{4m}$ which assigns two numbers to a decorated ideal triangle. For two decorated ideal triangulation $\tau$ and $\tau'$ the corresponding coordinate changes are rational.

For a decorated ideal triangulation $\tau$ and a number $q= \mathrm e ^{\pi \mathrm i \hbar} \in \mathbb C$, Kashaev introduced an algebra $\mathcal K_\tau^q$ which is the algebra over $\mathbb C$ defined by generators $Y_1^{\pm1}$, $Z_1^{\pm1}$, \dots, $Y_{2m}^{\pm1}$, $Z_{2m}^{\pm1}$ associated to ideal triangles of $\tau$ and by relations
\begin{align*}
Y_iY_j&=Y_jY_i, \\
Z_iZ_j&=Z_jZ_i, \\
Y_iZ_j&=Z_jY_i\ \ \mbox{if}\ \ i\neq j, \\
Z_iY_i&=q^2Y_iZ_i.
\end{align*}
Let $\widehat{\mathcal K}_{\tau}^q $ be the fraction division algebra of $\mathcal K_\tau^q$.

As one moves from one decorated ideal triangulation $\tau$ to another $\tau'$, Kashaev \cite{Kas1} introduce coordinate change isomorphisms from $\widehat{\mathcal K}_{\tau'}^q$ to $\widehat{\mathcal K}_\tau^q$. Analog to the construction of quantum Teichm\"uller space, there is an algebra $\widehat{\mathcal K}_S^q$ associated to a surface which is independent of decorated ideal triangulations.

In this paper, we will show Kashaev's construction of coordinate change isomorphisms are not unique. In fact, for two arbitrary complex numbers $a, b$, there are coordinate change isomorphisms from $\widehat{\mathcal K}_{\tau'}^q$ to $\widehat{\mathcal K}_\tau^q$. And Kashaev's construction is the special case of $a=q^{-1},b=q.$ For this generalized coordinate change isomorphisms, we also obtain a well-defined noncommutative
algebra $\widehat{\mathcal K}_S^q (a,b)$ associated to the surface $S$ which is called the \emph{generalized Kashaev algebra}. This is stated in Theorem \ref{thm:main}.

\subsection{Relationship between quantum Teichm\"uller space and Kashaev algebra}

To understand the relationship between the quantum Teichm\"uller space and Kashaev algebra, we need to first understand the relationship between shear coordinates and Kashaev coordinates. Fix a decorated ideal triangulation, the space of Kashaev coordinates is a fiber bundle on a subset in the enhanced Teichm\"uller space whose fiber is an affine space modeled on the homology group $H_1(S,\mathbb{R})$. This is proved in Theorem \ref{thm:exact}.

The relationship between the shear coordinates and Kashaev coordinates induces a natural relationship between the quantum Teichm\"uller space $\widehat{\mathcal T}_S^q$ and the generalized Kashaev algebra $\widehat{\mathcal K}_S^q (a,b)$. We show that there is a homomorphism from the quotient algebra $\widehat{\mathcal T}_S^q/(q^{-2m-\sum_{i<j}\sigma^\lambda_{ij}}X_1X_2...X_{3m})$ to $\widehat{\mathcal K}_S^q (a,b)$ if and only if $a=q^{-2}$ and $b=q^3$. This is proved in Corollary \ref{cor:homo} and Theorem \ref{thm:homo}. The result explains why we need to look for new construction of coordinate change isomorphisms other than Kashaev's construction.

\subsection{Open questions}

Hua Bai \cite{Bai} proved that the construction of quantum Teichm\"uller space $\widehat{\mathcal T}_S^q$ is essentially unique. The uniqueness of the algebra $\widehat{\mathcal K}_S^q (a,b)$ should be an interesting problem.

In \cite{BonLiu, Liu2, Bon2}, it is shown that quantum Teichm\"uller space $\widehat{\mathcal T}_S^q$ has a rich representation theory and an invariant of hyperbolic 3-manifolds is constructed. The representation theory of the algebra $\widehat{\mathcal K}_S^q (a,b)$ should be investigated.

The motivation of our work is to understand the relationship between Chekhov-Fock and Kashaev's operator-theoretical quantization. It is important to find a relationship between the two operator algebras involving the quantum dilogarithm function.

\section{Decorated ideal triangulations}

Let $S$ be an oriented surface of genus $g$ with $p\geq 1$ punctures and negative Euler characteristic, i.e., $m=2g-2+p>0.$  Any ideal triangulation has $2m$ ideal triangles and $3m$ edges.

A decorated ideal triangulation $\tau$ of $S$ is introduced by Kashaev \cite{Kas1} as an ideal triangulation such that the ideal triangles are numerated as $\{\tau_1,\tau_2,...,\tau_{2m}\}$ and there is a mark (a star symbol) at a corner of each ideal triangle. Denote by $\triangle(S)$ the set of isotopy classes of decorated ideal triangulations of surface $S$.

The set $\triangle(S)$ admits a natural action of the group
$\mathfrak{S}_{2m}$ of permutations of $2m$ elements, acting by
permuting the indices of the ideal triangles of $\tau$. Namely $\tau'
= \alpha(\tau)$ for $\alpha \in\mathfrak{S}_{2m}$ if its $i$--th
ideal triangle $\tau'_i$ is equal to $\tau_{\alpha(i)}$.

\begin{figure}[ht!]
\labellist\small\hair 2pt
\pinlabel $(1)$ at 205 200
\pinlabel $\longrightarrow$ at 205 81
\pinlabel $\varphi_{ij}$ at 205 99
\pinlabel $\tau_i$ at 50 85
\pinlabel $\tau_j$ at 107 85
\pinlabel $*$ at 13 85
\pinlabel $*$ at 144 85
\pinlabel $\tau'_i$ at 329 107
\pinlabel $\tau'_j$ at 329 48
\pinlabel $*$ at 329 144
\pinlabel $*$ at 329 13
\pinlabel $(2)$ at 205 0
\pinlabel $\longrightarrow$ at 205 273
\pinlabel $\rho_i$ at 205 291
\pinlabel $\tau_i$ at 80 264
\pinlabel $\tau'_i$ at 329 264
\pinlabel $*$ at 79 301
\pinlabel $*$ at 260 240

\endlabellist
\centering
\includegraphics[scale=0.5]{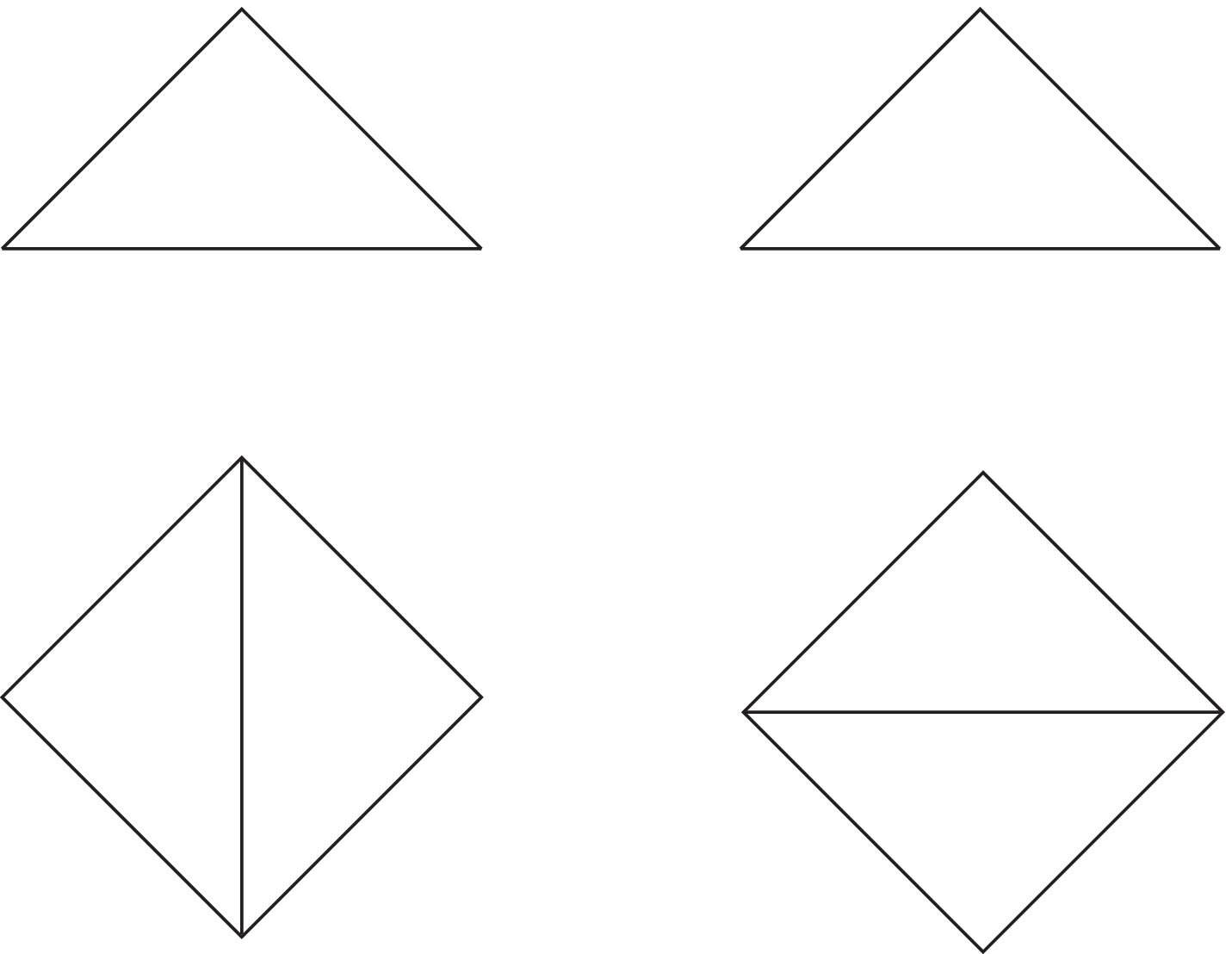}
\caption{}
\label{fig:maps}
\end{figure}

Another important transformation of $\triangle(S)$ is provided by
the \emph {diagonal exchange} $\varphi_{ij}:
\triangle(S)\rightarrow \triangle(S)$ defined as follows. Suppose that two ideal triangles $\tau_i, \tau_j$ share an edge $e$ such that the marked corners are opposite to the edge $e.$ Then $\varphi_{ij}(\tau)$ is obtained by rotating the interior of the union $\tau_i\cup \tau_j$ $90^\circ$ in the clockwise order, as illustrated in Figure \ref{fig:maps}(2).

The last one of transformations of $\triangle(S)$ is the \emph {mark rotation} $\rho_i: \triangle(S)\rightarrow \triangle(S)$. $\rho_i(\tau)$ is obtained by relocating the mark of the ideal triangle $\tau_i$ from one corner to the next corner in the counterclockwise order, as illustrated in Figure \ref{fig:maps}(1).

\begin{figure}[ht!]
\labellist\small\hair 2pt
\pinlabel $\tau_i$ at 33 81
\pinlabel $\tau_j$ at 81 42
\pinlabel $\tau_k$ at 127 81
\pinlabel $*$ at 7 92
\pinlabel $*$ at 42 7
\pinlabel $*$ at 133 17

\endlabellist
\centering
\includegraphics[scale=0.5]{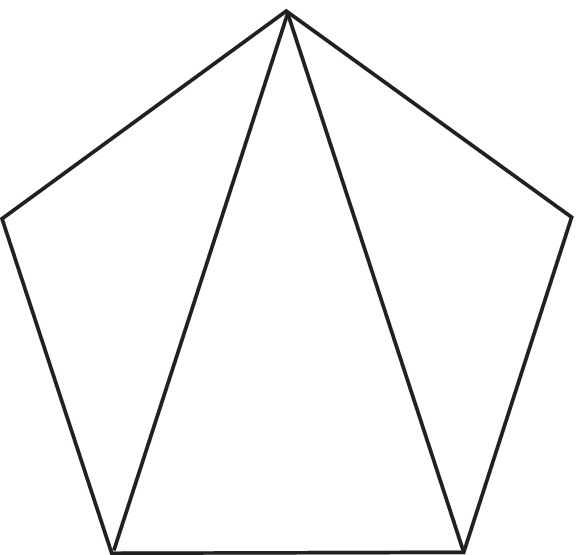}
\caption{}
\label{fig:pentagon}
\end{figure}

\begin{lemma}\label{thm:graph}
The reindexings, diagonal exchanges and mark rotations satisfy the following
relations:
\begin{enumerate}

\item $(\alpha\beta)(\tau) = \alpha(\beta(\tau))$ for every $\alpha$, $\beta\in\mathfrak{S}_{2m}$;

\item $\varphi_{ij}\circ \varphi_{ij}=\alpha_{i\leftrightarrow j},$ where $\alpha_{i\leftrightarrow j}$ denotes the transposition exchanging $i$ and $j$;

\item $\alpha \circ \varphi_{ij} = \varphi_{\alpha(i)\alpha(j)}\circ \alpha$ for every $\alpha\in\mathfrak{S}_{2m}$;

\item $\varphi_{ij}\circ\varphi_{kl}(\tau)=\varphi_{kl}\circ \varphi_{ij}(\tau)$, for $\{i,j\}\neq\{k,l\}$;

\item If three triangles $\tau_i,\tau_j,\tau_k$ of an ideal
triangulation $\tau \in \triangle(S)$ form a pentagon and their marked corners are in the location as in Figure \ref{fig:pentagon}, then the Pentagon Relation holds:
$$\omega_{jk}\circ\omega_{ik}\circ\omega_{ij}(\tau)=\omega_{ij}\circ\omega_{jk}(\tau),$$
where $\omega_{\mu\nu}=\rho_\mu\circ\varphi_{\mu\nu}\circ\rho_\nu$;

\item $\rho_i \circ \rho_i \circ \rho_i=\mathrm{Id}$;

\item $\rho_i \circ \rho_j= \rho_j \circ \rho_i$;

\item $\alpha \circ \rho_i= \rho_i \circ \alpha$ for for every $\alpha\in\mathfrak{S}_{2m}$.

\end{enumerate}

\end{lemma}

The lemma can be proved by drawing graphs.

\noindent\textbf{Remark.} Lemma \ref{thm:graph} is essential contained in Kashaev \cite{Kas2} where $\omega_{ij}$ is used as the diagonal exchange.
\medskip

The following two results about decorated ideal triangulations can be easily proved using Penner's result about ideal triangulations \cite{Pen}.

\begin{theorem}\label{thm:Penner1} Given two decorated ideal triangulations
$\tau,\tau'\in \triangle(S)$, there exists a finite sequence
of decorated ideal triangulations $\tau=\tau_{(0)}$, $\tau_{(1)}$,
\dots, $\tau_{(n)}=\tau'$ such that each $\tau_{(k+1)}$
is obtained from $\tau_{(k)}$ by a diagonal exchange or by a mark rotation or by a
reindexing of its ideal triangles.
\end{theorem}

\begin{theorem}\label{thm:Penner2} Given two decorated ideal triangulations
$\tau,\tau'\in \triangle(S)$ and given two sequences
$\tau=\tau_{(0)}$, $\tau_{(1)}$, \ldots, $\tau_{(n)}
=\tau'$ and $\tau=\overline{\tau}_{(0)}$, $\overline{\tau}_{(1)}$, \ldots,
$\overline{\tau}_{(\overline{n})} =\tau'$ of diagonal exchanges, mark rotations and reindexings
connecting them as in Theorem \ref{thm:Penner1}, these two
sequences can be related to each other by successive applications
of the following moves and of their inverses. These moves correspond to the relation in Lemma \ref{thm:graph}.
\begin{enumerate}

\item Replace \dots,
$\tau_{(k)}$, $\beta(\tau_{(k)})$, $\alpha\circ \beta(\tau_{(k)})$, \dots

\ \ \ by \dots, $\tau_{(k)}$, $(\alpha \beta)(\tau_{(k)})$, \dots
where $\alpha$, $\beta\in \mathfrak{S}_n$.

\item Replace \dots,
$\tau_{(k)}$, $\varphi_{ij}(\tau_{(k)}) $, $\varphi_{ij}\circ \varphi_{ij}(\tau_{(k)}) $ \dots

\ \ \ by \dots, $ \tau_{(k)}$, $\alpha_{i\leftrightarrow j}(\tau_{(k)})$, \dots .

\item Replace \dots,
$\tau_{(k)}$, $\varphi_{ij}(\tau_{(k)}) $, $\alpha \circ \varphi_{ij}(\tau_{(k)})$, \dots

\ \ \ by \dots, $\tau_{(k)}$, $\alpha(\tau_{(k)})$, $\varphi_{\alpha(i)\alpha(j)}\circ \alpha(\tau_{(k)})$, \dots
where $\alpha\in \mathfrak{S}_n$.

\item Replace \dots, $\tau_{(k)}$, $\varphi_{kl}(\tau_{(k)})$, $\varphi_{ij}\circ \varphi_{kl}(\tau_{(k)})$, \dots

\ \ \ by \dots, $ \tau_{(k)}$, $\varphi_{ij}(\tau_{(k)})$, $\varphi_{kl}\circ \varphi_{ij}(\tau_{(k)})$, \dots  where $\{i,j\}\neq\{k,l\}$.

\item Replace \dots, $\tau_{(k)}$, $\omega_{ij}(\tau_{(k)})$, $\omega_{ik} \circ\omega_{ij}(\tau_{(k)})$,
$\omega_{jk} \circ\omega_{ik} \circ\omega_{ij}(\tau_{(k)})$, \dots,

\ \ \ by \dots, $\tau_{(k)}$, $\omega_{jk}(\tau_{(k)})$, $\omega_{ij} \circ\omega_{jk}(\tau_{(k)})$, \dots
where $\omega_{\mu\nu}=\rho_\mu\circ\varphi_{\mu\nu}\circ\rho_\nu$.

\item Replace \dots, $\tau_{(k)}$,  $\rho_i(\tau_{(k)})$,  $\rho_i\circ \rho_i(\tau_{(k)})$, $\tau_{(k)}$ \dots

\ \ \ by  \dots, $\tau_{(k)}$, \dots.

\item Replace \dots, $\tau_{(k)}$, $\rho_i(\tau_{(k)})$, $\rho_j(\tau_{(k)})$, \dots

\ \ \ by \dots, $\tau_{(k)}$, $\rho_j(\tau_{(k)})$, $\rho_i(\tau_{(k)})$, \dots.

\item Replace \dots, $\tau_{(k)}$, $\rho_i(\tau_{(k)})$, $\alpha\circ \rho_i(\tau_{(k)})$, \dots

\ \ \ by \dots, $\tau_{(k)}$, $\alpha(\tau_{(k)})$, $\rho_i\circ \alpha(\tau_{(k)})$, \dots.

\end{enumerate}
\end{theorem}

\section{Generalized Kashaev algebra}

For a decorated ideal triangulation $\tau$ of a punctured surface $S$, Kashaev \cite{Kas1} associated each ideal triangle $\tau_i$ two numbers $\ln y_i, \ln z_i$. A Kashaev coordinate is a vector $(\ln y_1, \ln z_1,...,\ln y_{2m}, \ln z_{2m})\in \mathbb{R}^{4m}.$

Denote by $(y_1,z_1,...,y_{2m},z_{2m})$ the exponential Kashaev coordinate for the decorated ideal triangulation $\tau$. Denote by $(y'_1,z'_1,...,y'_{2m},z'_{2m})$ the exponential Kashaev coordinate for the decorated ideal triangulation $\tau'$.

\begin{definition}[Kashaev \cite{Kas1}]\label{def:coor-change}
Suppose that a decorated ideal triangulation $\tau'$ is obtained from another one $\tau$ by reindexing the ideal triangles, i.e., $\tau'=\alpha(\tau)$ for some $\alpha\in \mathfrak{S}_{2m},$ then we define $(y'_i, z'_i)=(y_{\alpha(i)}, z_{\alpha(i)})$ for any $i=1,...,2m.$

Suppose that a decorated ideal triangulation $\tau'$ is obtained from another one $\tau$ by a mark rotation, i.e., $\tau'=\rho_i(\tau)$ for some $i,$ then we define $(y'_j, z'_j)=(y_j,z_j)$ for any $j\neq i$ while
$$(y'_i, z'_i)=(\frac{z_i}{y_i}, \frac1{y_i}).$$

Suppose a decorated ideal triangulation $\tau'$ is obtained from another one $\tau$ by a diagonal exchange, i.e., $\tau'=\varphi_{ij}(\tau)$ for some $i,j,$ then we define $(y'_k, z'_k)=(y_k,z_k)$ for any $k\notin \{i,j\}$ while
$$(y'_i,z'_i,y'_j,z'_j)=
(\frac{z_j}{y_iy_j+z_iz_j},\frac{y_i}{y_iy_j+z_iz_j},\frac{z_i}{y_iy_j+z_iz_j},\frac{y_j}{y_iy_j+z_iz_j}).$$
\end{definition}

\noindent\textbf{Remark.} Kashaev \cite{Kas1} considered $\omega_{ij}$ instead of $\varphi_{ij}$.
\medskip

There is a natural relationship between Kashaev coordinates and Penner coordinates which is established in \cite{Kas1}. For exposition, see Teschner \cite{Te}. In Appendix, we include the main feature of this topic. Especially, the changes of Kashaev coordinate in Definition \ref{def:coor-change} are compatible with the changes of Penner coordinates.

For a decorated ideal triangulation $\tau$ of a punctured surface $S$, Kashaev \cite{Kas1} introduced an algebra $\mathcal K^q_{\tau}$ on $\mathbb C$ generated by $Y_1^{\pm},Z_1^{\pm},Y_2^{\pm},Z_2^{\pm},...,Y_{2m}^{\pm},Z_{2m}^{\pm},$ with $Y_i^{\pm},Z_i^{\pm}$ associated to an ideal triangle $\tau_i,$ and by the relations:
\begin{equation}\label{fml:kalgebra}
\begin{split}
Y_iY_j&=Y_jY_i, \\
Z_iZ_j&=Z_jZ_i, \\
Y_iZ_j&=Z_jY_i\ \ \mbox{if}\ \ i\neq j, \\
Z_iY_i&=q^2Y_iZ_i
\end{split}
\end{equation}

\noindent\textbf{Remark.} Kashaev's original definition is $Y_iZ_i=q^2Z_iY_i.$ We adopt our convention to make it compatible with the quantum Teichm\"uller space \cite{Liu1}. Kashaev's parameter $q$ is our $q^{-1}$.
\medskip

The algebra $\widehat{\mathcal K}^q_{\tau}$ is the fraction division algebra of $\mathcal K^q_{\tau}$ which consists of all the factors $FG^{-1}$ with $F,G\in \mathcal K^q_{\tau}$ and $Q\neq 0,$ and two such fractions
$F_1G_1^{-1}$ and $F_2G_2^{-1}$ are identified if there exists
$S_1$, $S_2\in \mathcal K^q_{\tau} - \{0\}$ such that $P_1S_1
= P_2S_2$ and $Q_1S_1 = Q_2S_2$.

In particular, when $q=1,$ $\mathcal K^q_{\tau}$ and $\widehat{\mathcal K}^q_{\tau}$ respectively coincide with the Laurent polynomial algebra $\mathbb C [Y_1^{\pm},Z_1^{\pm},...,Y_{2m}^{\pm},Z_{2m}^{\pm}]$ and the rational
fraction algebra $\mathbb C (Y_1,Z_1,...,Y_{2m},Z_{2m})$. The general
$\mathcal K^q_{\tau}$ and $\widehat{\mathcal K}^q_{\tau}$ can be considered as deformations of $\mathcal K^1_{\tau}$ and $\widehat{\mathcal K}^1_{\tau}$.

The algebra $\widehat{\mathcal K}^q_{\tau}$ depends on the decorated ideal triangulation $\tau$. We introduce algebra isomorphisms in the following.

\begin{definition}\label{def:iso} For any numbers $a,b\in \mathbb C.$

Suppose that a decorated ideal triangulation $\tau'$ is obtained from another one $\tau$ by reindexing the ideal triangles, i.e., $\tau'=\alpha(\tau)$ for some $\alpha\in \mathfrak{S}_{2m},$ then we define a map $\widehat{\alpha}:\widehat{\mathcal K}^q_{\tau'}\to \widehat{\mathcal K}^q_{\tau}$ by indicating the image of generators and extend it to the whole algebra:
\begin{align*}
\widehat{\alpha}(Y'_i) &= Y_{\alpha(i)}, \ \ \mbox{for any}\ \ i=1,...,2m,\\
\widehat{\alpha}(Z'_i) &= Z_{\alpha(i)}, \ \ \mbox{for any}\ \ i=1,...,2m.
\end{align*}

Suppose that a decorated ideal triangulation $\tau'$ is obtained from another one $\tau$ by a mark rotation, i.e., $\tau'=\rho_i(\tau)$ for some $i,$ then we define a map $\widehat{\rho}_i:\widehat{\mathcal K}^q_{\tau'}\to \widehat{\mathcal K}^q_{\tau}$ by indicating the image of generators and extend it to the whole algebra:
\begin{align*}
\widehat{\rho}_i(Y'_j) &= Y_j, \ \ \mbox{if}\ \ j\neq i, \\
\widehat{\rho}_i(Z'_j) &= Z_j, \ \ \mbox{if}\ \ j\neq i, \\
\widehat{\rho}_i(Y'_i) &= aY_i^{-1}Z_i, \\
\widehat{\rho}_i(Z'_i) &= Y_i^{-1}.
\end{align*}

Suppose a decorated ideal triangulation $\tau'$ is obtained from another one $\tau$ by a diagonal exchange, i.e., $\tau'=\varphi_{ij}(\tau)$ for some $i,j,$ then we define a map $\widehat{\varphi}_{ij}:\widehat{\mathcal K}^q_{\tau'}\to \widehat{\mathcal K}^q_{\tau}$ by indicating the image of generators and extend it to the whole algebra:
\begin{align*}
\widehat{\varphi}_{ij}(Y_i')&= \ (b Y_iY_j+Z_iZ_j)^{-1}Z_j, \\
\widehat{\varphi}_{ij}(Z_i')&=b(b Y_iY_j+Z_iZ_j)^{-1}Y_i, \\
\widehat{\varphi}_{ij}(Y_j')&= \ (b Y_iY_j+Z_iZ_j)^{-1}Z_i, \\
\widehat{\varphi}_{ij}(Z_j')&=b(b Y_iY_j+Z_iZ_j)^{-1}Y_j.
\end{align*}

\end{definition}

\noindent\textbf{Remark.} From the definition, when $a=b=1,$ we get the coordinate change formula in Definition \ref{def:coor-change}.

\noindent\textbf{Remark.} Kashaev \cite{Kas1} considered a special case of these maps when $a=q^{-1}, b=q.$

\begin{prop}\label{thm:relation} The maps $\widehat{\alpha}, \widehat{\rho}_i$ and $\widehat{\varphi}_{ij}$ satisfy the following relations which correspond to the relations in Lemma \ref{thm:graph}:

\begin{enumerate}

\item $\widehat{\alpha\beta}= \widehat{\alpha}\circ\widehat{\beta}$ for every $\alpha$, $\beta\in\mathfrak{S}_{2m}$;

\item $\widehat{\varphi}_{ij}\circ\widehat{\varphi}_{ij}=\widehat{\alpha}_{i\leftrightarrow j}$;

\item $ \widehat{\alpha} \circ \widehat{\varphi}_{ij}  =
 \widehat{\varphi}_{\alpha(i)\alpha(j)} \circ \widehat{\alpha} $ for every
$\alpha\in\mathfrak{S}_{2m}$;

\item $\widehat{\varphi}_{ij}\circ\widehat{\varphi}_{kl}=\widehat{\varphi}_{kl}\circ \widehat{\varphi}_{ij}$ for $\{i,j\}\neq\{k,l\}$;

\item If three triangles $\tau_i,\tau_j,\tau_k$ of an ideal
triangulation $\tau \in \triangle(S)$ form a pentagon and their marked corners are in the location as in Figure \ref{fig:pentagon}, then the Pentagon Relation holds:
$$\widehat{\omega}_{jk}\circ\widehat{\omega}_{ik}\circ\widehat{\omega}_{ij}=
\widehat{\omega}_{ij}\circ\widehat{\omega}_{jk},$$
where $\widehat{\omega}_{\mu\nu}=\widehat{\rho}_\mu\circ\widehat{\varphi}_{\mu\nu}\circ\widehat{\rho}_\nu$;

\item $\widehat{\rho}_i\circ\widehat{\rho}_i\circ\widehat{\rho}_i=\mathrm{Id}$;

\item $\widehat{\rho}_i\circ\widehat{\rho}_j=\widehat{\rho}_j\circ\widehat{\rho}_i$;

\item $\widehat{\alpha} \circ \widehat{\rho}_i= \widehat{\rho}_i \circ \widehat{\alpha}$ for every
$\alpha\in\mathfrak{S}_{2m}$.

\end{enumerate}
\end{prop}

\begin{proof} (1),(3),(4),(7),(8) are obvious.

(6) can be proved by using definition of $\widehat{\rho}_i$ easily. In fact, we assume that
$$\tau \xleftarrow{\rho_i}  \tau'' \xleftarrow{\rho_i}  \tau'\xleftarrow{\rho_i}  \tau.$$ Then we have
$$\widehat{\mathcal{K}}^q_{\tau} \xrightarrow{\widehat{\rho}_i}  \widehat{\mathcal{K}}^q_{\tau''} \xrightarrow{\widehat{\rho}_i} \widehat{\mathcal{K}}^q_{\tau'} \xrightarrow{\widehat{\rho}_i}  \widehat{\mathcal{K}}^q_{\tau}.$$
To show (6) is true, we need to show that $\widehat{\rho}_i\circ\widehat{\rho}_i\circ\widehat{\rho}_i$ sends each generator of $\widehat{\mathcal{K}}^q_{\tau}$ to itself. This is true for $Y_j,Z_j,j\neq i.$ We only need to take care of $Y_i,Z_i.$
For example, we check that $\widehat{\rho}_i\circ\widehat{\rho}_i\circ\widehat{\rho}_i(Y_i)=Y_i.$
In fact,
\begin{align*}
Y_i\xrightarrow{\widehat{\rho}_i}\ \ & aY''^{-1}_iZ''_i \\
\xrightarrow{\widehat{\rho}_i}\ \ &  a\widehat{\rho}_i(Y''^{-1}_i)\widehat{\rho}_i(Z''_i)=a(aY'^{-1}_iZ'_i)^{-1}Y'^{-1}_i=Z'^{-1}_i\\
\xrightarrow{\widehat{\rho}_i}\ \ & \widehat{\rho}_i(Z'_i)^{-1}=Y_i.
\end{align*}

To prove (2), we assume that
$$\tau \xleftarrow{\alpha_{i\leftrightarrow j}}  \tau'' \xleftarrow{\varphi_{ij}}  \tau'\xleftarrow{\varphi_{ij}}  \tau.$$ Then we have
$$\widehat{\mathcal{K}}^q_{\tau} \xrightarrow{\widehat{\alpha}_{i\leftrightarrow j}}  \widehat{\mathcal{K}}^q_{\tau''} \xrightarrow{\widehat{\varphi}_{ij}} \widehat{\mathcal{K}}^q_{\tau'} \xrightarrow{\widehat{\varphi}_{ij}}  \widehat{\mathcal{K}}^q_{\tau}.$$
To show that $\widehat{\varphi}_{ij}\circ\widehat{\varphi}_{ij}\circ \widehat{\alpha}_{i\leftrightarrow j}=\mathrm{Id},$ we need to show that it sends every generator of $\widehat{\mathcal{K}}^q_{\tau}$ to itself. This is true for $Y_k,Z_k,k\notin \{i,j\}.$ We only need to take care of $Y_i,Z_i.Y_j,Z_j.$
For example, we check it for $Y_i.$ In fact,
\begin{align*}
Y_i \xrightarrow{\widehat{\alpha}_{i\leftrightarrow j}} \ \ & Y''_j \\
\xrightarrow{\widehat{\varphi}_{ij}}\ \ &  (b Y'_iY'_j+Z'_iZ'_j)^{-1}Z'_i,\\
\xrightarrow{\widehat{\varphi}_{ij}}\ \ &
[b(b Y_iY_j+Z_iZ_j)^{-1}Z_j(b Y_iY_j+Z_iZ_j)^{-1}Z_i+\\
& b^2(b Y_iY_j+Z_iZ_j)^{-1}Y_i(b Y_iY_j+Z_iZ_j)^{-1}Y_j]^{-1}\widehat{\varphi}_{ij}(Z'_i)\\
=\ \ & [b(b Y_iY_j+Z_iZ_j)^{-1}(b q^2Y_iY_j+Z_iZ_j)^{-1}Z_jZ_i+\\
&b^2(b Y_iY_j+Z_iZ_j)^{-1}(bq^2 Y_iY_j+Z_iZ_j)^{-1}q^2Y_iY_j]^{-1}\widehat{\varphi}_{ij}(Z'_i)\\
=\ \ &[b(b Y_iY_j+Z_iZ_j)^{-1}(bq^2 Y_iY_j+Z_iZ_j)^{-1}(Z_jZ_i+bq^2Y_iY_j)]^{-1}\widehat{\varphi}_{ij}(Z'_i)\\
=\ \ & b^{-1}(b Y_iY_j+Z_iZ_j)\widehat{\varphi}_{ij}(Z'_i)\\
=\ \ & b^{-1}(b Y_iY_j+Z_iZ_j)b(b Y_iY_j+Z_iZ_j)^{-1}Y_i\\
=\ \ & Y_i.
\end{align*}

To prove (4) the Pentagon Relation, we need more work. As stated in Lemma \ref{thm:graph}, the Pentagon Relation for the decorated ideal triangulation is
\begin{align*}
& \omega_{jk}\circ\omega_{ik}\circ\omega_{ij}=\omega_{ij}\circ\omega_{jk}\\
\Longleftrightarrow\ \ \ & \rho_j\circ\varphi_{jk}\circ\rho_k\circ\rho_i\circ\varphi_{ik}\circ\rho_k
\circ\rho_i\circ\varphi_{ij}\circ\rho_j=\rho_i\circ\varphi_{ij}\circ\rho_j\circ\rho_j\circ\varphi_{jk}\circ\rho_k\\
\Longleftrightarrow\ \ \ &
\end{align*}
\begin{align*}
\rho_j\circ\varphi_{jk}\circ\rho_k\circ\rho_i\circ\varphi_{ik}\circ\rho_k
\circ\rho_i\circ\varphi_{ij}\circ\rho_j \circ \rho_k^2\circ \alpha_{j\leftrightarrow k}\circ \varphi_{jk}
\circ \rho_j \circ \alpha_{i\leftrightarrow j} \circ \varphi_{ij}\circ \rho_i^2 = \mathrm{Id},
\end{align*}
since $\rho_i^{-1}=\rho_i^2$ and $\varphi_{ij}^{-1}=\alpha_{i\leftrightarrow j} \circ \varphi_{ij}.$
Assume that
\begin{align*}
\tau        \xleftarrow{\rho_j}
\tau^{(17)} \xleftarrow{\varphi_{jk}}
\tau^{(16)} \xleftarrow{\rho_k}
\tau^{(15)} \xleftarrow{\rho_i}
\tau^{(14)} \xleftarrow{\varphi_{ik}}
\tau^{(13)} \\
            \xleftarrow{\rho_k}
\tau^{(12)} \xleftarrow{\rho_i}
\tau^{(11)} \xleftarrow{\varphi_{ij}}
\tau^{(10)} \xleftarrow{\rho_j}
\tau^{(9)}  \xleftarrow{\rho_k}
\tau^{(8)} \xleftarrow{\rho_k}
\tau^{(7)} \\
           \xleftarrow{\alpha_{j\leftrightarrow k}}
\tau^{(6)} \xleftarrow{\varphi_{jk}}
\tau^{(5)} \xleftarrow{\rho_j}
\tau^{(4)} \xleftarrow{\alpha_{i\leftrightarrow j}}
\tau^{(3)} \xleftarrow{\varphi_{ij}}
\tau^{(2)} \xleftarrow{\rho_i}
\tau^{(1)} \xleftarrow{\rho_i}
\tau
\end{align*}
Then we have
\begin{align*}
\widehat{\mathcal{K}}^q_{\tau}       \xrightarrow{\widehat{\rho}_j}
\widehat{\mathcal{K}}^q_{\tau^{(17)}} \xrightarrow{\widehat{\varphi}_{jk}}
\widehat{\mathcal{K}}^q_{\tau^{(16)}} \xrightarrow{\widehat{\rho}_k}
\widehat{\mathcal{K}}^q_{\tau^{(15)}} \xrightarrow{\widehat{\rho}_i}
\widehat{\mathcal{K}}^q_{\tau^{(14)}} \xrightarrow{\widehat{\varphi}_{ik}}
\widehat{\mathcal{K}}^q_{\tau^{(13)}} \\
            \xrightarrow{\widehat{\rho}_k}
\widehat{\mathcal{K}}^q_{\tau^{(12)}} \xrightarrow{\widehat{\rho}_i}
\widehat{\mathcal{K}}^q_{\tau^{(11)}} \xrightarrow{\widehat{\varphi}_{ij}}
\widehat{\mathcal{K}}^q_{\tau^{(10)}} \xrightarrow{\widehat{\rho}_j}
\widehat{\mathcal{K}}^q_{\tau^{(9)}}  \xrightarrow{\widehat{\rho}_k}
\widehat{\mathcal{K}}^q_{\tau^{(8)}} \xrightarrow{\widehat{\rho}_k}
\widehat{\mathcal{K}}^q_{\tau^{(7)}} \\
           \xrightarrow{\widehat{\alpha}_{j\leftrightarrow k}}
\widehat{\mathcal{K}}^q_{\tau^{(6)}} \xrightarrow{\widehat{\varphi}_{jk}}
\widehat{\mathcal{K}}^q_{\tau^{(5)}} \xrightarrow{\widehat{\rho}_j}
\widehat{\mathcal{K}}^q_{\tau^{(4)}} \xrightarrow{\widehat{\alpha}_{i\leftrightarrow j}}
\widehat{\mathcal{K}}^q_{\tau^{(3)}} \xrightarrow{\widehat{\varphi}_{ij}}
\widehat{\mathcal{K}}^q_{\tau^{(2)}} \xrightarrow{\widehat{\rho}_i}
\widehat{\mathcal{K}}^q_{\tau^{(1)}} \xrightarrow{\widehat{\rho}_i}
\widehat{\mathcal{K}}^q_{\tau}.
\end{align*}
To verify the Pentagon Relation, we need to show that the composition of maps sends every generator of $\widehat{\mathcal{K}}^q_{\tau}$ to itself. This is true for $Y_l,Z_l,l\notin\{i,j,k\}$. We only need to take care of $Y_i,Z_i,Y_j,Z_j,Y_k,Z_k.$ For example, we verify it holds for $Y_i.$
For simplicity the notation, in the following calculation, we do not indicate the upper index of generators. For example the second $Y_i$ should be $Y_i^{(17)}.$

\begin{align*}
Y_i  \xrightarrow{\widehat{\rho}_j} \ \ & Y_i \\
\xrightarrow{\widehat{\varphi}_{jk}} \ \ & Y_i \\
\xrightarrow{\widehat{\rho}_k} \ \ & Y_i \\
\xrightarrow{\widehat{\rho}_i} \ \ & aY^{-1}_iZ_i\\
\xrightarrow{\widehat{\varphi}_{ik}} \ \ & a[(bY_iY_k+Z_iZ_k)^{-1}Z_k]^{-1}b(bY_iY_k+Z_iZ_k)^{-1}Y_i=abZ^{-1}_kY_i\\
\xrightarrow{\widehat{\rho}_k} \ \ & abY_kY_i\\
\xrightarrow{\widehat{\rho}_i} \ \ & a^2bY_kY^{-1}_iZ_i\\
\xrightarrow{\widehat{\varphi}_{ij}} \ \ & a^2bY_k[(bY_iY_j+Z_iZ_j)^{-1}Z_j]^{-1}b(bY_iY_j+Z_iZ_j)^{-1}Y_i=a^2b^2Y_kZ_j^{-1}Y_i \\
\xrightarrow{\widehat{\rho}_j} \ \ & a^2b^2Y_kY_jY_i \\
\xrightarrow{\widehat{\rho}_k} \ \ & a^3b^2Y^{-1}_kZ_kY_jY_i \\
\xrightarrow{\widehat{\rho}_k} \ \ & a^2b^2Z^{-1}_kY_jY_i \\
\xrightarrow{\widehat{\alpha}_{j\leftrightarrow k}} \ \ & a^2b^2Z^{-1}_jY_kY_i \\
\xrightarrow{\widehat{\varphi}_{jk}} \ \ & a^2b^2[b(bY_jY_k+Z_jZ_k)^{-1}Y_j]^{-1}(bY_jY_k+Z_jZ_k)^{-1}Z_jY_i=a^2bY^{-1}_jZ_jY_i\\
\xrightarrow{\widehat{\rho}_j} \ \ & abZ^{-1}_jY_i \\
\xrightarrow{\widehat{\alpha}_{i\leftrightarrow j}} \ \ & abZ^{-1}_iY_j \\
\xrightarrow{\widehat{\varphi}_{ij}} \ \ & ab[b(bY_iY_j+Z_iZ_j)^{-1}Y_i]^{-1}(bY_iY_j+Z_iZ_j)^{-1}Z_i=aY^{-1}_iZ_i\\
\xrightarrow{\widehat{\rho}_i} \ \ & Z^{-1}_i\\
\xrightarrow{\widehat{\rho}_i} \ \ & Y_i.
\end{align*}

\end{proof}

\begin{prop}\label{thm:iso} If a decorated ideal triangulation $\tau'$ is obtained from another one $\tau$ by an operation $\pi,$ where $\pi=\alpha$ for some $\alpha\in \mathfrak{S}_{2m},$ or $\pi=\rho_i$ for some $i$, or $\pi=\varphi_{ij}$ for some $i,j$, then $\widehat{\pi}: \widehat{\mathcal K}^q_{\tau'}\to \widehat{\mathcal K}^q_{\tau}$ as in Definition \ref{def:iso} is an isomorphism between the two algebras.
\end{prop}

\begin{proof} If $\pi=\alpha$ for some $\alpha\in \mathfrak{S}_{2m},$ it is obvious that $\widehat{\pi}$ is an isomorphism.

If $\pi=\rho_i$ for some $i$, we need to check that $\widehat{\pi}$ is a homomorphism, i.e., it preserve the algebraic relations (\ref{fml:kalgebra}). The first three are obvious. It is enough to check the last one. Since $Z'_iY'_i=q^2Y'_iZ'_i,$ we need to show that $\widehat{\pi}(Z'_iY'_i)=q^2\widehat{\pi}(Y'_iZ'_i).$ We verify this by showing

\begin{align*} & &\widehat{\pi}(Z'_iY'_i) &=q^2\widehat{\pi}(Y'_iZ'_i)\\
&\Longleftrightarrow  &\widehat{\pi}(Z'_i)\widehat{\pi}(Y'_i) &=q^2\widehat{\pi}(Y'_i)\widehat{\pi}(Z'_i) \\
&\Longleftrightarrow  &Y_i^{-1}aY_i^{-1}Z_i &=q^2aY_i^{-1}Z_iY_i^{-1}
                           \ \ \ \ \ \mbox{by Definition \ref{def:iso}}\\
&\Longleftrightarrow  &Y_i^{-1}Z_i &=q^2 Z_iY_i^{-1} \\
&\Longleftrightarrow  &Z_iY_i &=q^2 Y_iZ_i.
\end{align*}
This is true.

To show that $\widehat{\rho}_i$ is an isomorphism, it is enough to find its inverse. In fact, by Proposition \ref{thm:relation}(6), we see
$\widehat{\rho}_i^{-1}=\widehat{\rho}_i\circ \widehat{\rho}_i.$

If $\pi=\varphi_{ij}$ for some $i,j$, we need to check that $\widehat{\pi}$ is a homomorphism, i.e., it preserve the algebraic relations (\ref{fml:kalgebra}).

\textbf{Case 1}: For $\{k,l\}\neq \{i,j\}$, since $\widehat{\pi}(Y'_k,Z'_k,Y'_l,Z'_l)=(Y_k,Z_k,Y_l,Z_l),$ therefore $\widehat{\pi}$ preserves the relation of $Y'_k,Z'_k,Y'_l,Z'_l$.

\textbf{Case 2}: For $k\notin \{i,j\}$, we consider $Y'_k,Z'_k$ and $Y'_i,Z'_i$. Now

\begin{align*}
\widehat{\varphi}_{ij}(Y_k')&=Y_k,\\
\widehat{\varphi}_{ij}(Z_k')&=Z_k,\\
\widehat{\varphi}_{ij}(Y_i')&= \ (b Y_iY_j+Z_iZ_j)^{-1}Z_j, \\
\widehat{\varphi}_{ij}(Z_i')&=b(b Y_iY_j+Z_iZ_j)^{-1}Y_i.
\end{align*}

Since $Y_k'Y_i'=Y_i'Y_k',$ we need to check that $\widehat{\pi}(Y_k'Y_i')=\widehat{\pi}(Y_i'Y_k').$ This is true.

Since $Z_k'Z_i'=Z_i'Z_k',$ we need to check that $\widehat{\pi}(Z_k'Z_i')=\widehat{\pi}(Z_i'Z_k').$ This is true.

Since $Z_k'Y_k'=q^2Y_k'Z_k',$ we need to check that $\widehat{\pi}(Z_k'Y_k')=q^2\widehat{\pi}(Y_k'Z_k').$ This is true.

Since $Z_i'Y_i'=q^2Y_i'Z_i',$ we need to check that $\widehat{\pi}(Z_i'Y_i')=q^2\widehat{\pi}(Y_i'Z_i')$ which is equivalent to
\begin{multline*}
b(b Y_iY_j+Z_iZ_j)^{-1}Y_i  (b Y_iY_j+Z_iZ_j)^{-1}Z_j\\
=q^2(b Y_iY_j+Z_iZ_j)^{-1}Z_j b(b Y_iY_j+Z_iZ_j)^{-1}Y_i
\end{multline*}
\begin{align*}
&\Longleftrightarrow & Y_i(b Y_iY_j+Z_iZ_j)^{-1}Z_j =  q^2  Z_j(b Y_iY_j+Z_iZ_j)^{-1}Y_i\\
&\Longleftrightarrow & Z_jY_i^{-1} (b Y_iY_j+Z_iZ_j) =q^2 (b Y_iY_j+Z_iZ_j)Y_i Z_j^{-1}.
\end{align*}
This is true since
\begin{align*}
\mbox{the left hand side} &=Z_j(b Y_iY_j+q^2Z_iZ_j)Y_i^{-1} \\
&= (bq^2 Y_iY_j+q^2Z_iZ_j)Z_jY_i^{-1}\\
&= q^2(b Y_iY_j+Z_iZ_j)Z_jY_i^{-1}\\
&= \mbox{the right hand side}.
\end{align*}

\textbf{Case 3}: We consider $Y'_i,Z'_i$ and $Y'_j,Z'_j$.

Since $Y_i'Y_j'=Y_j'Y_i',$ we need to check that $\widehat{\pi}(Y_i'Y_j')=\widehat{\pi}(Y_j'Y_i')$ which is equivalent to
\begin{multline*}
(b Y_iY_j+Z_iZ_j)^{-1}Z_j  (b Y_iY_j+Z_iZ_j)^{-1}Z_i\\
=(b Y_iY_j+Z_iZ_j)^{-1}Z_j  (b Y_iY_j+Z_iZ_j)^{-1}Y_i
\end{multline*}
\begin{align*}
&\Longleftrightarrow & Z_j(b Y_iY_j+Z_iZ_j)^{-1}Z_i =  Z_i(b Y_iY_j+Z_iZ_j)^{-1}Z_j\\
&\Longleftrightarrow & Z_iZ_j^{-1} (b Y_iY_j+Z_iZ_j) =(b Y_iY_j+Z_iZ_j)Z_j^{-1}Z_i.
\end{align*}
This is true since
\begin{align*}
\mbox{the left hand side} &= Z_i(bq^{-2} Y_iY_j+Z_iZ_j)Z_j^{-1} \\
&= (bq^{-2}q^2 Y_iY_j+Z_iZ_j)Z_iZ_j^{-1}\\
&= \mbox{the right hand side}.
\end{align*}

The similar calculation is used to check that $\widehat{\pi}(Z_i'Z_j')=\widehat{\pi}(Z_j'Z_i')$.

Since $Y_i'Z_j'=Z_j'Y_i',$ we need to check that $\widehat{\pi}(Y_i'Z_j')=\widehat{\pi}(Z_j'Y_i')$ which is equivalent to
\begin{multline*}
(b Y_iY_j+Z_iZ_j)^{-1}Z_j  b(b Y_iY_j+Z_iZ_j)^{-1}Y_j\\
=b(b Y_iY_j+Z_iZ_j)^{-1}Y_j (b Y_iY_j+Z_iZ_j)^{-1}Z_j
\end{multline*}
\begin{align*}
&\Longleftrightarrow & Z_j(b Y_iY_j+Z_iZ_j)^{-1}Y_j =  Y_j(b Y_iY_j+Z_iZ_j)^{-1}Z_j\\
&\Longleftrightarrow & Y_jZ_j^{-1} (b Y_iY_j+Z_iZ_j) =(b Y_iY_j+Z_iZ_j)Z_j^{-1}Y_j.
\end{align*}
This is true since
\begin{align*}
\mbox{the left hand side} &= Y_j(bq^{-2} Y_iY_j+Z_iZ_j)Z_j^{-1} \\
&= (bq^{-2}Y_iY_j+q^{-2}Z_iZ_j)Y_jZ_j^{-1}\\
&= (bq^{-2}Y_iY_j+q^{-2}Z_iZ_j)q^2Z_j^{-1}Y_j\\
&= \mbox{the right hand side}.
\end{align*}

The similar calculation is used to check that $\widehat{\pi}(Y_j'Z_i')=\widehat{\pi}(Z_i'Y_j')$.

For $Z_i'Y_i'=q^2Y_i'Z_i'$ and $Z_j'Y_j'=q^2Y_j'Z_j'$, we have done in Case 2.

To show that $\widehat{\varphi}_{ij}$ is an isomorphism, it is enough to find its inverse. In fact, by Proposition \ref{thm:relation}(2), we see $\widehat{\varphi}_{ij}^{-1}=\widehat{\alpha}_{i\leftrightarrow j}\circ \widehat{\varphi}_{ij},$ where $\alpha_{i\leftrightarrow j}$ denotes the transposition exchanging $i$ and $j$.
\end{proof}

\begin{theorem}\label{thm:main}
For two arbitrary complex numbers $a,b$, there is a unique family of algebra isomorphisms
$$\Psi_{\tau\tau'}^q(a,b):
\widehat{\mathcal{K}}^q_{\tau'} \rightarrow
\widehat{\mathcal{K}}^q_{\tau}$$ defined as $\tau$,
$\tau' \in \triangle(S)$ ranges over all pairs of decorated ideal
triangulations, such that:
\begin{enumerate}
\item $\Psi_{\tau\tau''}^q(a,b) = \Psi_{\tau\tau'}^q(a,b) \circ
\Psi_{\tau'\tau''}^q(a,b)$ for every $\tau$, $\tau'$,
$\tau''\in \triangle(S)$;

\item $\Psi_{\tau\tau'}^q(a,b)$ is the isomorphism of
Definition \ref{def:iso} when $\tau'$ is obtained
from $\tau$ by a reindexing or a mark rotation or a diagonal exchange.
\end{enumerate}
\end{theorem}

\begin{proof} Use Theorem \ref{thm:Penner1} to
connect $\tau$ to $\tau'$ by a sequence
$\tau=\tau_{(0)}$, $\tau_{(1)}$, \dots,
$\tau_{(n)}=\tau'$ where each $\tau_{(k+1)}$ is obtained
from $\tau_{(k)}$ by a reindexing or a mark rotation or a diagonal exchange,
and define $\Psi_{\tau\tau'}^q(a,b)$ as the
composition of the $\Psi_{\tau_{(k)}\tau_{(k+1)}}^q(a,b)$. Theorem \ref{thm:Penner2} and
Proposition \ref{thm:relation}
show that this $\Psi_{\tau_{(k)}\tau_{(k+1)}}^q(a,b)$ is
independent of the choice of the sequence of $\tau_{(k)}$.

The uniqueness immediately follows from Theorem \ref{thm:Penner1}.
\end{proof}

The \emph{generalized Kashaev algebra} $\widehat{\mathcal{K}}^q_S(a,b)$ associated to a surface $S$ is defined as the algebra
$$
\widehat{\mathcal{K}}^q_S(a,b)= \bigg(
\bigsqcup_{\tau \in \triangle(S)}
\widehat{\mathcal{K}}^q_\tau(a,b)\bigg)/\sim
$$
where the relation $\sim$ is defined by the property that, for
$X\in \widehat{\mathcal{K}}^q_\tau(a,b)$ and $X'\in
\widehat{\mathcal{K}}^q_{\tau'}(a,b)$,
$$
X \sim X' \Leftrightarrow X=\Psi^q_{\tau,\tau'}(a,b)(X').
$$

\section{Kashaev coordinates and shear coordinates}

To understand the quantization using shear coordinates and the quantization using Kashaev coordinates, we first need to understand the relationship between these two coordinates.

\subsection{Decorated ideal triangulations}
Given a decorated ideal triangulation $\tau,$ by forgetting the mark at each corner, we obtain an ideal triangulation $\lambda.$ We call $\lambda$ the \emph{underlying ideal triangulation} of $\tau$. let $\lambda_1,\lambda_2,...,\lambda_{3m}$ be the components of ideal triangulation $\lambda$. Denote by $\tau_1,..,\tau_{2m}$ the ideal triangles.

For the ideal triangulation $\lambda$, we may consider its dual graph. Each ideal triangle $\tau_\mu$ corresponds to a vertex $\tau_\mu^*$ of the dual graph. Denote by $\lambda_1^*,\lambda_2^*,...,\lambda_{3m}^*$ the dual edges. If an edge $\lambda_i$ bounds one side of the ideal triangles $\tau_\mu$ and one side of $\tau_\nu$, then the dual edge $\lambda_i^*$ connects the vertexes $\tau_\mu^*$ and $\tau_\nu^*$.

In a decorated ideal triangulation $\tau$, each ideal triangle $\tau_\mu$ (embedded or not) has three sides which correspond to the three half-edges incident to the vertex $\tau_\mu^*$ of the dual graph. The three sides are numerated by $0,1,2$ in the counterclockwise order such that the $0-$side is opposite to the marked corner.

\subsection{Space of Kashaev coordinates}
Let's recall that a Kashaev coordinate associated to a decorated ideal triangulation $\tau$ is a vector $(\ln y_1,\ln z_1,...,\ln y_{2m},\ln z_{2m})\in \mathbb{R}^{4m}$, where $\ln y_\mu$ and $\ln z_\mu$ are associated to the ideal triangle $\tau_\mu$. Denote by $\mathcal K_\tau$ the space of Kahsaev coordinates associated to $\tau$. We see that $\mathcal K_\tau=\mathbb{R}^{4m}.$

Given a vector $(\ln y_1,\ln z_1,...,\ln y_{2m},\ln z_{2m})\in \mathcal K_\tau$, we associate a number to each side of each ideal triangle as follows. For the ideal triangle $\tau_\mu$, we associate
\begin{enumerate}
\item[]$\ln h_\mu^0:=\ln y_\mu-\ln z_\mu$ to the 0-side;

\item[]$\ln h_\mu^1:=\ln z_\mu$ to the 1-side;

\item[]$\ln h_\mu^2:=-\ln y_\mu$ to the 2-side.
\end{enumerate}

Therefore $\ln h_\mu^0+\ln h_\mu^1+\ln h_\mu^2=0$. We can identify the space  $\mathcal K_\tau=\mathbb{R}^{4m}$ with a subspace of $\mathbb{R}^{6m}=\{(...,\ln h_\mu^0,\ln h_\mu^1,\ln h_\mu^2,...)\}$ satisfying $\ln h_\mu^0+\ln h_\mu^1+\ln h_\mu^2=0$ for each ideal triangle $\tau_\mu$.

\subsection{Exact sequence}
The enhanced Teichm\"uller space parametrized by shear coordinates is $\widetilde{\mathcal T}_\lambda=\mathbb{R}^{3m}=\{(\ln x_1,\ln x_2,...,\ln x_{3m})\}$, where $\ln x_i$ is the shear coordinate at edge $\lambda_i$. We define a map $f_1: \widetilde{\mathcal T}_\lambda \to \mathbb{R}$ by sending $(\ln x_1,\ln x_2,...,\ln x_{3m})$ to the sum of entries $\sum_{i=1}^{3m}\ln x_i.$

Suppose $\lambda$ is the underlying ideal triangulation of the decorated ideal triangulation $\tau$. We define a map $f_2: \mathcal K_\tau \to \widetilde{\mathcal T}_\lambda$ as a linear function by setting $$\ln x_i=\ln h_\mu^s+\ln h_\nu^t$$ whenever $\lambda_i$ bounds the $s-$side of $\tau_\mu$ and the $t-$side of $\tau_\nu$ ($\mu$ may equal $\nu$), where $s,t\in\{0,1,2\}$.

Another map $f_3: H_1(S,\mathbb{R}) \to \mathcal K_\tau$ is defined as follows. A homology class in $H_1(S,\mathbb{R})$ is represented by a linear combination of oriented dual edges: $\sum_{i=1}^{3m}c_i\lambda_i^*$. If the orientation of $\lambda_i^*$ is from the $s-$side of $\tau_\mu$ to the $t-$side of $\tau_\nu$, by setting $\ln h_\mu^s=-c_i$ and $\ln h_\nu^t=c_i$, we obtain a vector $(...,\ln h_\mu^0,\ln h_\mu^1,\ln h_\mu^2,...)\in \mathbb{R}^{6m}.$ The boundary map of chain complexes sends $\sum_{i=1}^{3m}c_i\lambda_i^*$ to a linear combination of vertexes. In this combination, the term involving the vertex $\tau_\mu^*$ is $(c_i\epsilon_i +c_j\epsilon_j+c_k\epsilon_k)\tau_\mu^*$ where $\lambda_i,\lambda_j,\lambda_k$ (two of them may coincide) bound three sides of $\tau_\mu$ and $\epsilon_t=-1$ if $\lambda_t^*$ starts at the side of $\tau_\mu$ bounded by $\lambda_t$ while $\epsilon_t=1$ if $\lambda_t^*$ ends at the side of $\tau_\mu$ bounded by $\lambda_t$. Therefore $$(c_i\epsilon_i +c_j\epsilon_j+c_k\epsilon_k)\tau_\mu^*=(\ln h_\mu^0+\ln h_\mu^1+\ln h_\mu^2)\tau_\mu^*.$$ Since the chain $\sum_{i=1}^{3m}c_i\lambda_i^*$ is a cycle, we must have $\ln h_\mu^0+\ln h_\mu^1+\ln h_\mu^2=0$. Therefore this vector $(...,\ln h_\mu^0,\ln h_\mu^1,\ln h_\mu^2,...)$ is in the subspace $\mathcal K_\tau$.

Combining the three maps, we obtain

\begin{theorem}\label{thm:exact} The following sequence is exact:
$$0\rightarrow H_1(S, \mathbb R)\xrightarrow{f_3} \mathcal K_\tau \xrightarrow{f_2} \widetilde{\mathcal T}_\lambda \xrightarrow{f_1} \mathbb R \rightarrow 0.$$
\end{theorem}

\begin{proof}
The map $f_3$ is injective. In fact, if the homology class represented by $\sum_{i=1}^{3m}c_i\lambda_i^*$ is mapped to the zero vector in $\mathcal K_\tau$, then, for each $i=1,...,3m,$ we have $|c_i|=|\ln h_\mu^s|=0,$ where $\lambda_i$ bounds the $s-$side of $\tau_\mu.$ Therefore it is a zero homology class. Thus the sequence is exact at $H_1(S, \mathbb R).$

Suppose $(...,\ln h_\mu^0,\ln h_\mu^1,\ln h_\mu^2,...)\in Im(f_3),$ i.e.,$$(...,\ln h_\mu^0,\ln h_\mu^1,\ln h_\mu^2,...)=f_3(\sum_{i=1}^{3m}c_i\lambda_i^*).$$ For any edge $\lambda_i$ bounds the $s-$side of $\tau_\mu$ and the $t-$side of $\tau_\nu$, we have $$\ln x_i= \ln h_\mu^s+ \ln h_\nu^t=\pm c_i\mp c_i=0.$$
Thus $(...,\ln h_\mu^i,\ln h_\mu^j,\ln h_\mu^k,...)\in Ker(f_2).$ Therefore $Im(f_3)\subseteq Ker(f_2).$

On the other hand, we claim $Im(f_3)\supseteq Ker(f_2)$. In fact, given a vector $(...,\ln h_\mu^0,\ln h_\mu^1,\ln h_\mu^2,...)\in Ker(f_2),$ we can reverse the process of the definition of $f_3$ to obtain a homology class in $H_1(S, \mathbb R)$. To be precise, since the vector is in the kernel of $f_2$, we have $\ln h_\mu^s+ \ln h_\nu^t=0$ for each edge $\lambda_i$ bounds the $s-$side of $\tau_\mu$ and the $t-$side of $\tau_\nu$. An orientation of $\lambda_i^*$ can be given as follows. 
\begin{enumerate}

\item[] When $\ln h_\mu^s > 0,$ the dual edge $\lambda_i^*$ runs from the $s-$side of $\tau_\mu$ to the $t-$side of $\tau_\nu$. 
    
\item[] When $\ln h_\mu^i< 0,$ the dual edge $\lambda_i^*$ runs from the $t-$side of $\tau_\nu$ to the $s-$side of $\tau_\mu$.
    
\item[] When $\ln h_\mu^s = 0,$ $\lambda_i^*$ is oriented in either way. 
\end{enumerate}

Consider the one dimensional chain $\sum_{i=1}^{3m}|\ln h_\mu^s|\lambda_i^*$, where $\lambda_i$ bounds the $s-$side of $\tau_\mu$. This chain turns out to be a cycle. In fact, the boundary map sends this chain to a zero dimensional chain in which the term involving the vertex $\tau_\mu^*$ is $$(|\ln h_\mu^0|\epsilon_0+|\ln h_\mu^1|\epsilon_1+|\ln h_\mu^2|\epsilon_2)\tau_\mu^*$$ where $\epsilon_s=\pm 1$ and $\epsilon_s=sign(\ln h_\mu^s)\cdot 1$ if $\ln h_\mu^s \neq 0$ for $s\in \{0,1,2\}$. Thus
$$(|\ln h_\mu^0|\epsilon_0+|\ln h_\mu^1|\epsilon_1+|\ln h_\mu^2|\epsilon_2)\tau_\mu^*=(\ln h_\mu^0+\ln h_\mu^1+\ln h_\mu^2)\tau_\mu^*=0.$$ This cycle defines a homology class.

The argument above shows $Im(f_3)= Ker(f_2)$, i.e., the sequence is exact at $\mathcal K_\tau$.

Now $\dim Ker(f_2)=\dim Im(f_3)=\dim H_1(S, \mathbb R)= 2g+p-1=m+1.$ Thus $\dim Im(f_2)=\dim \mathcal K_\tau - \dim Ker(f_2)=4m-(m+1)=3m-1.$ Since $$Ker(f_1)=\{(x_1,x_2,...,x_{3m})|\sum_{i=1}^{3m}x_i=0\}$$ is a subspace of dimension $3m-1$, then $Im(f_2)=Ker(f_1),$ i.e., the sequence is exact at $\widetilde{\mathcal T}_\lambda$.

It is easy to see that $f_1$ is onto. Therefore the sequence is exact at $\mathbb R$.

\end{proof}

\noindent\textbf{Remark.} From the theorem above, we see that $\mathcal K_\tau$ is a fiber bundle over the space $Ker(f_1)$ whose fiber is an affine space modeled on $H_1(S, \mathbb R).$ To be precise, given a vector $s\in Ker(f_1),$ let $v\in f_2^{-1}(s).$ Then $f_2^{-1}(s)=v+H_1(S, \mathbb R).$

\noindent\textbf{Remark.} There is an exact sequence relating space of Kashaev coordinates and decorated Teichm\"uller space. See Proposition \ref{thm:k-exact} in Appendix.
\medskip

\subsection{Relation of bivecotrs}

Consider the linear isomorphism
\begin{align}\label{fml:3-2}
M:  \mathcal K_\tau &\longrightarrow  \mathcal K_\tau\\
(\ln y_1,\ln z_1,...,\ln y_{2m},\ln z_{2m})&\longmapsto  (...,\ln h_\mu^0,\ln h_\mu^1,\ln h_\mu^2,...) .\notag
\end{align}

\begin{prop}\label{thm:bivector}
If $(\ln x_1,\ln x_2,...,\ln x_{3m})=f_2\circ M (\ln y_1,\ln z_1,...,\ln y_{2m},\ln z_{2m}),$ then
$$\sum_{i,j=1}^{3m}\sigma^\lambda_{ij}\frac{\partial}{\partial \ln x_i} \wedge \frac{\partial}{\partial \ln x_j}=
(f_2)_*\circ M_* (\sum_{\mu=1}^{2m}\frac{\partial}{\partial \ln y_\mu} \wedge \frac{\partial}{\partial \ln z_\mu}),$$
where $\sigma^\lambda_{ij}=a^\lambda_{ij}-a^\lambda_{ji}$ and $a^\lambda_{ij}$ is the number of corners of the ideal triangulation $\lambda$ which is delimited in the left by $\lambda_i$ and on the right by $\lambda_j$.
\end{prop}

\begin{proof} By definition (\ref{fml:3-2}), we have
$$M_*(\frac{\partial}{\partial \ln y_\mu} \wedge \frac{\partial}{\partial \ln z_\mu})=
\frac{\partial}{\partial \ln h_\mu^0} \wedge \frac{\partial}{\partial \ln h_\mu^1}+
\frac{\partial}{\partial \ln h_\mu^1} \wedge \frac{\partial}{\partial \ln h_\mu^2}+
\frac{\partial}{\partial \ln h_\mu^2} \wedge \frac{\partial}{\partial \ln h_\mu^1}.$$

Assume that the edges $\lambda_i,\lambda_j,\lambda_k$ (two of them may coincide) bound the $0-$side, $1-$side and $2-$side of the ideal triangle $\tau_\mu$ respectively.

By definition of map $f_2$, we have
\begin{align*}
(f_2)_*(\frac{\partial}{\partial \ln h_\mu^0} \wedge \frac{\partial}{\partial \ln h_\mu^1}+
\frac{\partial}{\partial \ln h_\mu^1} \wedge \frac{\partial}{\partial \ln h_\mu^2}+
\frac{\partial}{\partial \ln h_\mu^2} \wedge \frac{\partial}{\partial \ln h_\mu^0})\\
=\frac{\partial}{\partial \ln x_i} \wedge \frac{\partial}{\partial \ln x_j}+
\frac{\partial}{\partial \ln x_j} \wedge \frac{\partial}{\partial \ln x_k}+
\frac{\partial}{\partial \ln x_k} \wedge \frac{\partial}{\partial \ln x_i}.
\end{align*}
Therefore
\begin{align*}
&(f_2)_*\circ M_* (\sum_{\mu=1}^{2m}\frac{\partial}{\partial \ln y_\mu} \wedge \frac{\partial}{\partial \ln z_\mu})\\
=  &\sum_{\mu=1}^{2m}(\frac{\partial}{\partial \ln x_i} \wedge \frac{\partial}{\partial \ln x_j}+
\frac{\partial}{\partial \ln x_j} \wedge \frac{\partial}{\partial \ln x_k}+
\frac{\partial}{\partial \ln x_k} \wedge \frac{\partial}{\partial \ln x_i})\\
& \mbox{where $\lambda_i,\lambda_j,\lambda_k$ bound the $0-$side, $1-$side and $2-$side of $\tau_\mu$}\\
=  &\sum_{i,j=1}^{3m}\sigma^\lambda_{ij}\frac{\partial}{\partial \ln x_i} \wedge \frac{\partial}{\partial \ln x_j}.
\end{align*}

\end{proof}

\noindent\textbf{Remark.} There is a relationship between a differential 2-form in Kashaev coordinates and the Weil-Peterson 2-form in Penner coordinates. See Proposition \ref{thm:k-diff} in Appendix.
\medskip

\subsection{Compatibility of coordinate changes}

\begin{figure}[ht!]
\labellist\small\hair 2pt
\pinlabel $*$ at 11 129
\pinlabel $*$ at 129 11
\pinlabel $*$ at 307 11
\pinlabel $*$ at 429 129
\pinlabel $\tau_\mu$ at 37 96
\pinlabel $\tau_\mu$ at 400 96
\pinlabel $\tau_\nu$ at 98 42
\pinlabel $\tau_\nu$ at 340 42
\pinlabel $\tau_\eta$ at -17 114
\pinlabel $\tau_\eta$ at 286 114
\pinlabel $\tau_\zeta$ at 153 114
\pinlabel $\tau_\zeta$ at 452 114
\pinlabel $\lambda_i$ at 76 57
\pinlabel $\lambda_j$ at 66 159
\pinlabel $\lambda_m$ at -13 70
\pinlabel $\lambda_l$ at 66 -15
\pinlabel $\lambda_k$ at 155 70
\pinlabel $\lambda'_i$ at 381 78
\pinlabel $\lambda'_j$ at 365 159
\pinlabel $\lambda'_m$ at 285 70
\pinlabel $\lambda'_l$ at 365 -15
\pinlabel $\lambda'_k$ at 453 70
\pinlabel $\longrightarrow$ at 216 60
\pinlabel $\varphi_{\mu\nu}$ at 216 73

\endlabellist
\centering
\includegraphics[scale=0.5]{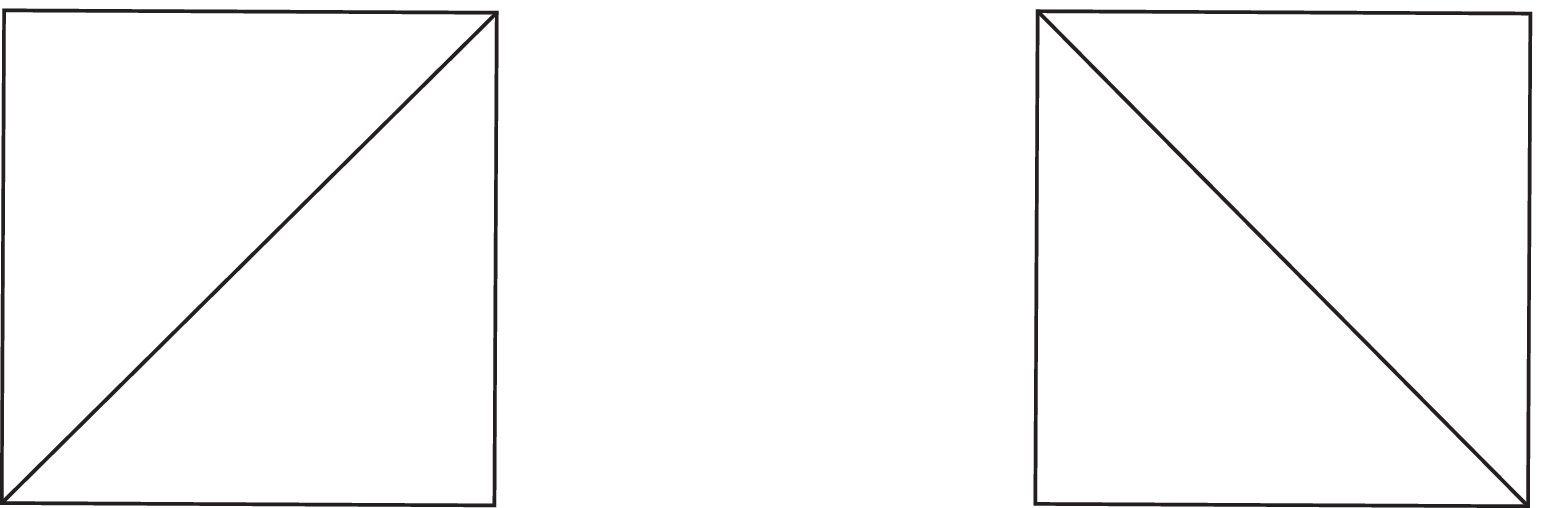}
\caption{}
\label{fig:eight}
\end{figure}

The coordinate change of shear coordinates are given as

\begin{prop}[\cite{Liu1} Proposition 3]
\label{prop:DiagExch} Suppose that the ideal triangulations
$\lambda$, $\lambda'$ are obtained from each other
by a diagonal exchange, namely that $\lambda' =
\Delta_i(\lambda)$. Label the edges of $\lambda$ involved in this
diagonal exchange as $\lambda_i$, $\lambda_j$, $\lambda_k$,
$\lambda_l$, $\lambda_m$ as in Figure~\ref {fig:eight}. If
$(x_1, x_2, \dots, x_n)$ and $(x_1', x_2', \dots,
x_n')$ are the exponential shear coordinates
associated $\lambda$ and $\lambda'$ of the same enhanced hyperbolic metric,
then $x_h'=x_h$ for every $h \not\in
\{i,j,k,l,m\}$, $x_i' = x_i^{-1}$ and:
\begin{description}
\item[Case 1] if the edges $\lambda_j$, $\lambda_k$, $\lambda_l$,
$\lambda_m$  are distinct, then
$$ x'_j = (1+x_i)x_j\quad x'_k  =
(1+x_i^{-1})^{-1}x_k \quad x'_l = (1+x_i)x_l \quad x'_m =
(1+x_i^{-1})^{-1}x_m;$$

\item[Case 2] if $\lambda_j$ is identified with $\lambda_k$, and
$\lambda_l$ is distinct from $\lambda_m$, then
$$
    x'_j  =  x_ix_j \quad x'_l  =
(1+x_i)x_l \quad x'_m  =  (1+x_i^{-1})^{-1}x_m;
$$

\item[Case 3] (the inverse of Case 2) if $\lambda_j$ is identified
with $\lambda_m$, and $\lambda_k$ is distinct from $\lambda_l$,
then
$$
    x'_j  =  x_ix_j
\quad x'_k  =  (1+x_i^{-1})^{-1}x_k \quad x'_l  = (1+x_i)x_l;
$$

\item[Case 4] if $\lambda_j$ is identified with $\lambda_l$, and
$\lambda_k$ is distinct from $\lambda_m$, then
$$
    x'_j  =  (1+x_i)^2x_j
\quad x'_k  =  (1+x_i^{-1})^{-1}x_k \quad x'_m  =
(1+x_i^{-1})^{-1}x_m
$$

\item[Case 5] (the inverse of Case 4) if $\lambda_k$ is identified
with $\lambda_m$, and $\lambda_j$ is distinct from $\lambda_l$,
then
$$ x'_j  =  (1+x_i)x_j \quad x'_k  =
(1+x_i^{-1})^{-2}x_k \quad x'_l = (1+x_i)x_l;
$$

\item[Case 6] if $\lambda_j$ is identified with $\lambda_k$, and
$\lambda_l$ is identified with $\lambda_m$ (in which case $S$ is a
$3$--times punctured sphere), then
$$
    x'_j  =  x_ix_j \quad x'_l  =
x_ix_l;
$$

\item[Case 7] (the inverse of Case 6) if $\lambda_j$ is identified
with $\lambda_m$, and $\lambda_k$ is identified with $\lambda_l$
(in which case $S$ is a $3$--times punctured sphere), then
$$
    x'_j  =  x_ix_j
\quad x'_k  =  x_ix_k;
$$

\item[Case 8] if $\lambda_j$ is identified with $\lambda_l$, and
$\lambda_k$ is identified with $\lambda_m$ (in which case $S$ is a
once punctured torus), then
$$ x'_j  =  (1+x_i)^2x_j \quad x'_k  =
(1+x_i^{-1})^{-2}x_k.
$$

\end{description}
\end{prop}

\begin{prop}\label{compatible} Suppose that the decorated ideal triangulations $\tau$ and $\tau'$ have the underlying ideal triangulations $\lambda$ and $\lambda'$ respectively. The following diagram is commutative:
$$
\begin{CD}
\widetilde{\mathcal{T}}_\lambda @<f_2<< \mathcal{K}_\tau\\
@VVV @VVV \\
\widetilde{\mathcal{T}}_{\lambda'} @<f_2<<  \mathcal{K}_{\tau'}
\end{CD}
$$
where the two vertical maps are corresponding coordinate changes. The coordinate changes of Kashaev coordinates are given in Definition \ref{def:coor-change}. The coordinate changes of shear coordinates are given in Proposition \ref{prop:DiagExch}.
\end{prop}

\begin{proof} For a reindexing, the conclusion is obvious. For a mark rotation, the conclusion is easily proved by definition. For diagonal exchange, we need to check the eight cases in Proposition \ref{prop:DiagExch}. For instance, we verify Case 4. As in Figure \ref{fig:eight}, through maps $f_2$ and $M$, we may identify
\begin{align*}
x_i &=\frac{y_\mu y_\nu}{z_\mu z_\nu}     &  x'_i &=\frac{y'_\mu y'_\nu}{z'_\mu z'_\nu}\\
x_j &=z_\mu z_\nu                         &  x'_j &= \frac1{y'_\mu y'_\nu} \\
x_k &=\frac{h^s_\zeta}{y_\nu}             &  x'_k &=z'_\mu h^s_\zeta \\
x_m &=\frac{h^t_\eta}{y_\mu}              &  x'_m &=z'_\nu h^t_\eta
\end{align*}
for some $s,t\in \{0,1,2\}.$

Then we have
\begin{align*}
\varphi_{\mu\nu}(x'_i)
=\varphi_{\mu\nu}(\frac{y'_\mu y'_\nu}{z'_\mu z'_\nu})
= \frac{z_\mu z_\nu}{y_\mu y_\nu}
=  x_i^{-1}.
\end{align*}

And
\begin{align*}
\varphi_{\mu\nu}(x'_j)
=\varphi_{\mu\nu}(\frac1{y'_\mu y'_\nu})
= \frac{(y_\mu y_\nu+z_\mu z_\nu)^2}{z_\mu z_\nu}
= (1+\frac{y_\mu y_\nu}{z_\mu z_\nu})^2 z_\mu z_\nu
= (1+x_i)^2x_j.
\end{align*}

And
\begin{align*}
\varphi_{\mu\nu}(x'_k)
=\varphi_{\mu\nu}(z'_\mu) h^s_\zeta
= \frac{y_\mu y_\nu}{y_\mu y_\nu+z_\mu z_\nu}\frac1{y_\mu}h^s_\zeta
= (1+x_i^{-1})^{-1}x_k.
\end{align*}

It is same for $x'_m$ due to the symmetry of $\mu,\nu.$

\end{proof}

\noindent\textbf{Remark.} The compatibility of coordinate changes of Kashaev coordinates and Penner coordinates is given in Appendix Propostion \ref{thm:k-com}.

\section{Relationship between quantum Teichm\"uller space and Kashaev algebra}

In this section, we establish a natural relationship between the quantum Teichm\"uller space $\widehat{\mathcal{T}}^q_S$ and the generalized Kashaev algebra $\widehat{\mathcal{K}}^q_S(a,b)$.

\subsection{Homomorphism}

For a decorated ideal triangulation $\tau$ of a punctured surface $S$, Kashaev \cite{Kas1} introduced an algebra $\mathcal K^q_{\tau}$ on $\mathbb C$ generated by $Y_1^{\pm},$ $Z_1^{\pm},$ $Y_2^{\pm},$ $Z_2^{\pm},...,Y_{2m}^{\pm},Z_{2m}^{\pm},$ with $Y_i^{\pm},Z_i^{\pm}$ associated to an ideal triangle $\tau_i,$ subject to the relations (\ref{fml:kalgebra}).

For a ideal triangle $\tau_\mu$, we associate three elements in $\mathcal{K}^q_{\tau}$ to the three sides of $\tau_\mu$ as follows:
\begin{enumerate}
\item[] $H^0_\mu:=Y_\mu Z^{-1}_\mu$ to the $0-$side;

\item[] $H^1_\mu:=Z_\mu$ to the $1-$side;

\item[] $H^2\mu:=Y^{-1}_\mu$ to the $2-$side.
\end{enumerate}

\begin{lemma}\label{thm:h} For any $s,t\in\{0,1,2\}$ and $\mu\in {1,2,...,3m},$
$$H^s_\mu H^t_\mu=q^{2\sigma_{st}}H^t_\mu H^s_\mu,$$
where $\sigma_{st}+\sigma_{ts}=0$ and $\sigma_{10}=\sigma_{02}=\sigma_{21}=1.$
\end{lemma}

\begin{proof} When $(s,t)=(1,0),$ we have $H^s_\mu=Z_\mu$ and $H^t_\mu=Y_\mu Z^{-1}_\mu.$ Thus
$$H^s_\mu H^t_\mu=Z_\mu Y_\mu Z^{-1}_\mu=q^2 Y_\mu Z^{-1}_\mu Z_\mu= q^2H^t_\mu H^s_\mu.$$

When $(s,t)=(0,2),$ we have $H^s_\mu=Y_\mu Z^{-1}_\mu$ and $H^t_\mu=Y^{-1}_\mu.$ Thus
$$H^s_\mu H^t_\mu=Y_\mu Z^{-1}_\mu Y^{-1}_\mu=q^2Y^{-1}_\mu Y_\mu Z^{-1}_\mu= q^2H^t_\mu H^s_\mu.$$

When $(s,t)=(2,1),$ we have $H^s_\mu=Y^{-1}_\mu$ and $H^t_\mu=Z_\mu.$ Thus
$$H^s_\mu H^t_\mu=Y^{-1}_\mu Z_\mu=q^2 Z_\mu Y^{-1}_\mu= q^2H^t_\mu H^s_\mu.$$
\end{proof}

Suppose $\lambda$ is the underlying ideal triangulation of $\tau$, the Chekhov-Fock algebra $\mathcal T_\lambda^q$ is the algebra over $\mathbb C$ defined by generators $X_1^{\pm1}$, $X_2^{\pm1}$, \dots, $X_n^{\pm1}$ associated to the components of $\lambda$ and by relations $X_iX_j=q^{2\sigma^\lambda_{ij}}X_jX_i$.

We define a map $F_\tau:\mathcal T_\lambda^q\to \mathcal K^q_{\tau}$ by indicating the image of the generators and extend it to the whole algebra. Suppose that the edge $\lambda_i$ bounds the $s-$side of $\tau_\mu$ and the $t-$side of $\tau_\nu$. We define

\begin{align}\label{def:X}
F_\tau(X_i)=q^{\delta_{\mu\nu}\sigma_{ts}}H^s_\mu H^t_\nu \in \mathcal{K}^q_\tau,
\end{align}

where $\sigma_{ts}$ is defined in Lemma \ref{thm:h} and $\delta_{\mu\nu}$ is the Kronecker delta, i.e., $\delta_{\mu\mu}=1$ and $\delta_{\mu\nu}=0$ if $\mu\neq \nu$. When $\mu=\nu,$ $X_i$ is well-defined, since $$q^{\sigma_{ts}}H^s_\mu H^t_\mu=q^{\sigma_{st}}H^t_\mu H^s_\mu$$ due to Lemma \ref{thm:h}.

This definition is natural since when $q=1$ we get the relationship between Kashaev coordinates and shear coordinates which is given by the map $M$ and $f_2$. In fact when $q=1$ then generators $Y_\mu, Z_\mu$ are commutative. They reduce to the geometric quantities $y_\mu,z_\mu$ associate to $\tau_\mu.$ $H^s_\mu$ and $X_i$ are reduced to $h^s_\mu$ and $x_i$.

\begin{lemma} The map $F_\tau:\mathcal T_\lambda^q\to \mathcal K^q_{\tau}$ is a homomorphism.
\end{lemma}

\begin{proof} It is enough to check $F_\tau(X_i)F_\tau(X_j)=q^{2\sigma^\lambda_{ij}}F_\tau(X_j)F_\tau(X_i)$ for any elements $X_i$ and $X_j$.
Assume the edge $\lambda_i$ bounds the $s-$side of $\tau_\mu$ and the $t-$side of $\tau_\nu$ while the edge $\lambda_j$ bounds the $k-$side of $\tau_\zeta$ and the $l-$side of $\tau_\eta$.

If $\{\mu,\nu\}\cap\{\zeta,\eta\}=\emptyset,$ then $F_\tau(X_i)$ commutes with $F_\tau(X_j)$. On the other hand, $\sigma^\lambda_{ij}=0.$ The statement holds.

If $(\mu,\nu,\zeta,\eta)=(\mu,\mu,\mu,\mu)$, then $X_i=X_j.$ The statement holds.

If $(\mu,\nu,\zeta,\eta)=(\mu,\mu,\mu,\eta), \mu\neq \eta,$ then $F_\tau(X_i)=q^{\sigma_{ts}}H^s_\mu H^t_\mu$ and $F_\tau(X_j)=H^k_\mu H^l_\eta.$ Thus by Lemma \ref{thm:h}, we have
$$F_\tau(X_i)F_\tau(X_j)=q^{2(\sigma_{tk}+\sigma_{sk})}F_\tau(X_j)F_\tau(X_i)=F_\tau(X_j)F_\tau(X_i)
=q^{2\sigma^\lambda_{ij}}F_\tau(X_j)F_\tau(X_i),$$
due to $\sigma^\lambda_{ij}=0$. 

If $(\mu,\nu,\zeta,\eta)=(\mu,\nu,\nu,\eta), \mu\neq \nu, \mu\neq \eta, \nu\neq \eta,$ then $F_\tau(X_i)=H^s_\mu H^t_\nu$ and $F_\tau(X_j)=H^k_\nu H^l_\eta.$ Thus by Lemma \ref{thm:h}, we have
$$F_\tau(X_i)F_\tau(X_j)=q^{2\sigma_{tk}}F_\tau(X_j)F_\tau(X_i)=q^{2\sigma^\lambda_{ij}}F_\tau(X_j)F_\tau(X_i).$$

If $(\mu,\nu,\zeta,\eta)=(\mu,\nu,\mu,\nu), \mu\neq \nu,$ then $F_\tau(X_i)=H^s_\mu H^t_\nu$ and $F_\tau(X_j)=H^k_\mu H^l_\nu.$ Thus by Lemma \ref{thm:h}, we have
$$F_\tau(X_i)F_\tau(X_j)=q^{2(\sigma_{sk}+\sigma_{tl})}F_\tau(X_j)F_\tau(X_i)
=q^{2\sigma^\lambda_{ij}}F_\tau(X_j)F_\tau(X_i).$$

\end{proof}

\subsection{Compatibility}

Chekhov-Fock algebra $\mathcal T_\lambda^q$ has a well-defined fraction division algebra $\widehat{\mathcal T}_\lambda^q$.
As one moves from one ideal triangulation $\lambda$ to another $\lambda'$, Chekhov and Fock \cite{Fo, FC, CF} (see also \cite{Liu1}) introduce coordinate change isomorphisms $\Phi_{\lambda\lambda'}^q: \widehat{\mathcal T}_{\lambda'} ^q \rightarrow \widehat{\mathcal T}_\lambda^q.$ 

\begin{prop}[\cite{Liu1} Proposition 5]\label{thm:exchange} Suppose that the ideal triangulations
$\lambda$, $\lambda'$ are obtained from each other
by a diagonal exchange, namely that $\lambda' =
\Delta_i(\lambda)$. Label the edges of $\lambda$ involved in this
diagonal exchange as $\lambda_i$, $\lambda_j$, $\lambda_k$,
$\lambda_l$, $\lambda_m$ as in Figure \ref {fig:eight}. Then
there is a unique algebra isomorphism
\begin{equation*}
\widehat{\Delta}_i: \widehat{\mathcal{T}}^q_{\lambda'}
\rightarrow \widehat{\mathcal{T}}^q_{\lambda}
\end{equation*}
such that
    $X_h' \mapsto X_h$
for every $h \not\in \{i,j,k,l,m\}$, $X_i' \mapsto X_i^{-1}$ and:
\begin{description}
\item[Case 1] if the edges $\lambda_j$, $\lambda_k$, $\lambda_l$,
$\lambda_m$  are distinct, then
\begin{align*}
X'_j &\mapsto (1+qX_i)X_j \qquad\, X'_k \mapsto
(1+qX_i^{-1})^{-1}X_k \\
    X'_l &\mapsto (1+qX_i)X_l
\qquad
    X'_m
\mapsto (1+qX_i^{-1})^{-1}X_m ;
\end{align*}

\item[Case 2] if $\lambda_j$ is identified with $\lambda_k$, and
$\lambda_l$ is distinct from $\lambda_m$, then
$$ X'_j  \mapsto   X_iX_j \quad X'_l
\mapsto (1+qX_i)X_l \quad X'_m \mapsto (1+qX_i^{-1})^{-1}X_m
$$

\item[Case 3] (the inverse of Case 2) if $\lambda_j$ is identified
with $\lambda_m$, and $\lambda_k$ is distinct from $\lambda_l$,
then
$$X'_j\mapsto   X_iX_j \quad X'_k
\mapsto (1+qX_i^{-1})^{-1}X_k \quad X'_l  \mapsto (1+qX_i)X_l
$$

\item[Case 4] if $\lambda_j$ is identified with $\lambda_l$, and
$\lambda_k$ is distinct from $\lambda_m$, then
\begin{gather*}
X'_j  \mapsto (1+qX_i) (1+q^3X_i) X_j \\ X'_k \mapsto
(1+qX_i^{-1})^{-1}X_k \quad X'_m \mapsto (1+qX_i^{-1})^{-1}X_m
\end{gather*}

\item[Case 5] (the inverse of Case 4) if $\lambda_k$ is identified
with $\lambda_m$, and $\lambda_j$ is distinct from $\lambda_l$,
then
\begin{gather*}
    X'_j  \mapsto  (1+qX_i)X_j
\quad X'_l
    \mapsto  (1+qX_i)X_l \\
    X'_k \mapsto
(1+qX_i^{-1})^{-1}(1+q^3X_i^{-1})^{-1}X_k
\end{gather*}

\item[Case 6] if $\lambda_j$ is identified with $\lambda_k$, and
$\lambda_l$ is identified with $\lambda_m$ (in which case $S$ is a
$3$--times punctured sphere), then
$$
    X'_j  \mapsto  X_iX_j \quad X'_l  \mapsto
X_iX_l;
$$

\item[Case 7] (the inverse of Case 6) if $\lambda_j$ is identified
with $\lambda_m$, and $\lambda_k$ is identified with $\lambda_l$
(in which case $S$ is a $3$--times punctured sphere), then
$$
    X'_j  \mapsto  X_iX_j
\quad X'_k  \mapsto  X_iX_k;
$$

\item[Case 8] if $\lambda_j$ is identified with $\lambda_l$, and
$\lambda_k$ is identified with $\lambda_m$ (in which case $S$ is a
once punctured torus), then
\begin{gather*}
    X'_j  \mapsto
(1+qX_i)(1+q^3X_i)X_j \\
    X'_k
\mapsto (1+qX_i^{-1})^{-1}(1+q^3X_i^{-1})^{-1}X_k
\end{gather*}

\end{description}

\end{prop}

Recall that $\widehat{\mathcal K}^q_{\tau}$ is the fraction division algebra of $\mathcal K^q_{\tau}$. The algebraic isomorphism between $\widehat{\mathcal K}^q_{\tau}$ and $\widehat{\mathcal K}^q_{\tau'}$ is defined in Definition \ref{def:iso}.

\begin{lemma}\label{thm:a}
Suppose that a decorated ideal triangulation $\tau'$ is obtained from $\tau$ by a mark rotation $\rho_\mu$ for some $\mu\in\{1,2,...,2m\}$. Let $\lambda$ be the common underlying ideal triangulation of $\tau$ and $\tau'$. The following diagram is commutative if and only if $a=q^{-2}$.
$$
\begin{CD}
\widehat{\mathcal T}_\lambda^q @>F_\tau>> \widehat{\mathcal K}^q_{\tau}\\
@A\mathrm{Id}AA @AA\widehat{\rho}_\mu A \\
\widehat{\mathcal T}_\lambda^q @>F_{\tau'}>>  \widehat{\mathcal K}^q_{\tau'}
\end{CD}
$$

\end{lemma}

\begin{proof} It is enough to check $F_\tau(X_i)=\widehat{\rho}_\mu\circ F_{\tau'}(X_i)$ holds for any generator $X_i$. 

If $\lambda_i$ does not bound a side of the ideal triangle $\tau_\mu,$ then $F_\tau(X_i)=\widehat{\rho}_\mu\circ F_{\tau'}(X_i)$ holds automatically.

Suppose $\lambda_i$ bounds the $s-$side of $\tau_\mu$ and the $t-$side of $\tau_\nu$. If $\mu\neq \nu$, then $\lambda_i$ bounds the $(s+2)$(modulo 3)-side of $\tau'_\mu$ and the $t-$side of $\tau'_\nu$. Then $F_\tau(X_i)=H^s_\mu H^t_\nu$ and $F_{\tau'}(X_i)=H'^{s+2}_\mu H^t_\nu$. To show $F_\tau(X_i)=\widehat{\rho}_\mu\circ F_{\tau'}(X_i)$ is enough to show that $H^s_\mu=\widehat{\rho}_\mu(H'^{s+2}_\mu)$.

If $\mu= \nu,$ then then $\lambda_i$ bounds the $(s+2)$(modulo 3)-side of $\tau'_\mu$ and the $(t+2)$(modulo 3)-side of $\tau'_\nu$. Then $F_\tau(X_i)=q^{\sigma_{ts}}H^s_\mu H^t_\mu$ and $F_{\tau'}(X_i)=q^{\sigma_{ts}}H'^{s+2}_\mu H'^{t+2}_\mu$. To show $F_\tau(X_i)=\widehat{\rho}_\mu\circ F_{\tau'}(X_i)$ is enough to show that $H^s_\mu=\widehat{\rho}_\mu(H'^{s+2}_\mu)$ for $s\in\{0,1,2\}$, since  $t\in\{0,1,2\}$.

When $s=0,$ we have $H^s_\mu=Y_\mu Z^{-1}_\mu$ and $H'^{s+2}_\mu=Y'^{-1}_\mu.$
Now
\begin{align*}
& &H^s_\mu&=\widehat{\rho}_\mu(H'^{s+2}_\mu)\\
&\Longleftrightarrow&Y_\mu Z^{-1}_\mu&=\widehat{\rho}_\mu(Y'_\mu)^{-1})\\
&\Longleftrightarrow&Y_\mu Z^{-1}_\mu&=a^{-1}Z^{-1}_\mu Y_\mu\\
&\Longleftrightarrow&Z_\mu Y_\mu &=a^{-1} Y_\mu Z_\mu\\
&\Longleftrightarrow&a&=q^{-2}.
\end{align*}

When $s=1$, we have $H^s_\mu=Z_\mu $ and $H'^{s+2}_\mu=Y'_\mu Z'^{-1}_\mu.$
Now
\begin{align*}
& &H^s_\mu&=\widehat{\rho}_\mu(H'^{s+2}_\mu)\\
&\Longleftrightarrow&Z_\mu&=\widehat{\rho}_\mu(Y'_\mu Z'^{-1}_\mu)\\
&\Longleftrightarrow&Z_\mu&=aY_\mu^{-1}Z_\mu Y_\mu \\
&\Longleftrightarrow&Z_\mu&=aq^2 Z_\mu\\
&\Longleftrightarrow&a&=q^{-2}.
\end{align*}

When $s=2$, we have $H^s_\mu=Y_\mu^{-1}$ and $H'^{s+2}_\mu=Z'_\mu.$
Now
\begin{align*}
& &H^i_\mu&=\widehat{\rho}_\mu(H'^i_\mu)\\
&\Longleftrightarrow&Y_\mu^{-1}&=\widehat{\rho}_\mu(Z'_\mu)\\
&\Longleftrightarrow&Y_\mu^{-1}&=Y_\mu^{-1}.
\end{align*}
This holds automatically.

\end{proof}

\begin{lemma}\label{thm:b}
Suppose that a decorated ideal triangulation $\tau'$ is obtained from $\tau$ by a diagonal exchange $\varphi_{\mu\nu}.$ Let $\lambda$ and $\lambda'$ be the underlying ideal triangulation of $\tau$ and $\tau'$ respectively. Then $\lambda'$ is obtained $\lambda$ by a diagonal exchange with respect to the edge $\lambda_i$ which is the common edge of $\tau_\mu$ and $\tau_\nu$.
The following diagram is commutative if and only if $b=q^{3}$.
$$
\begin{CD}
\widehat{\mathcal T}_\lambda^q @>F_\tau>> \widehat{\mathcal K}^q_{\tau}\\
@A\widehat{\Delta}_i AA @AA\widehat{\varphi}_{\mu\nu} A \\
\widehat{\mathcal T}_{\lambda'}^q @>F_{\tau'}>>  \widehat{\mathcal K}^q_{\tau'}
\end{CD}
$$
\end{lemma}

\begin{proof}

\begin{figure}[ht!]
\labellist\small\hair 2pt
\pinlabel $*$ at 11 89
\pinlabel $*$ at 124 11
\pinlabel $*$ at 279 9
\pinlabel $*$ at 397 89
\pinlabel $\tau_\mu$ at 52 59
\pinlabel $\tau_\mu$ at 353 59
\pinlabel $\tau_\nu$ at 98 15
\pinlabel $\tau_\nu$ at 303 15
\pinlabel $\tau_\eta$ at 83 127
\pinlabel $\tau_\eta$ at 323 127
\pinlabel $\lambda_i$ at 98 39
\pinlabel $\lambda_i$ at 303 41
\pinlabel $\lambda_j$ at 81 89
\pinlabel $\lambda_j$ at 323 89
\pinlabel $\longrightarrow$ at 201 63
\pinlabel $\varphi_{\mu\nu}$ at 201 73

\endlabellist
\centering
\includegraphics[scale=0.5]{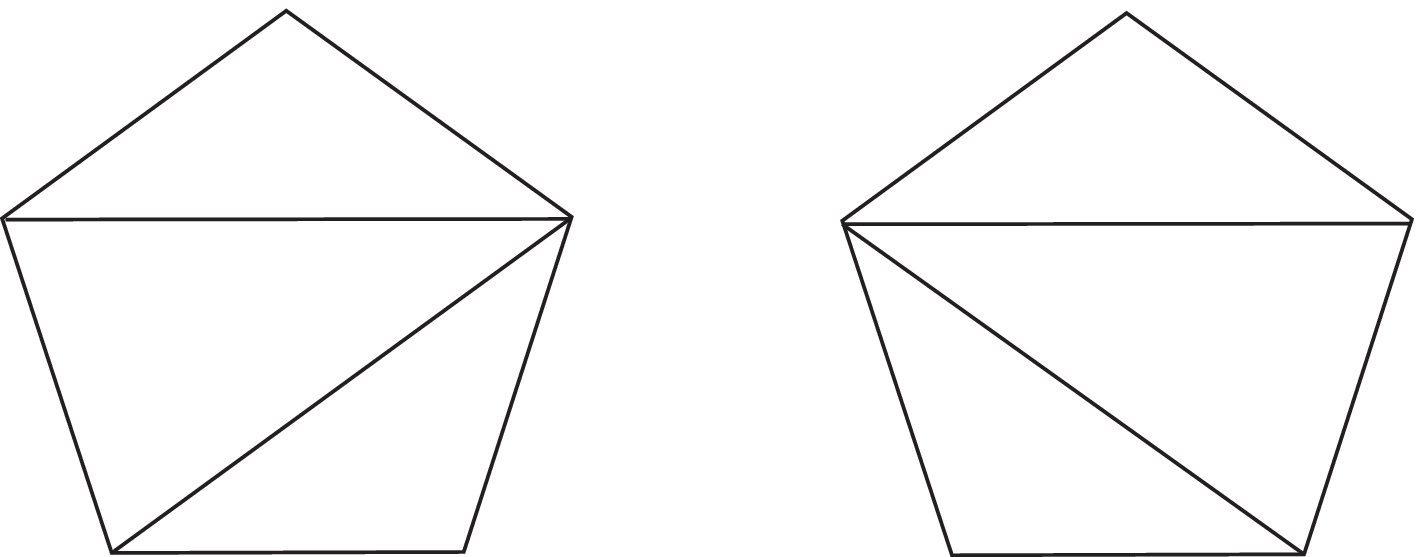}
\caption{}
\label{fig:b}
\end{figure}

First we show that $b=q^3$ is necessary. As in Figure \ref{fig:b}, $\mu,\nu,\eta$ are different. We have $F_\tau(X_i)=Y_\mu Z^{-1}_\mu Y_\nu Z^{-1}_\nu$ and $F_\tau(X_j)=H^s_\eta Z_\mu,$ for some $s\in \{0,1,2\}$.
And $F_{\tau'}(X_j)=H^j_\eta Y'^{-1}_\mu.$ That the diagram is commutative implies 
\begin{align*}
& & \widehat{\varphi}_{\mu\nu}\circ F_{\tau'}(X_j)&= F_\tau \circ \widehat{\Delta}_i(X'_j)\\
&\Longleftrightarrow &\widehat{\varphi}_{\mu\nu}(H^j_\eta Y'^{-1}_\mu)&=F_\tau((1+qX_i)X_j)\\
&\Longleftrightarrow & H^s_\eta [(bY_\mu Y_\nu+Z_\mu Z_\nu)^{-1}Z_\nu]^{-1}
                     &=(1+qY_\mu Z^{-1}_\mu Y_\nu Z^{-1}_\nu)H^s_\eta Z_\mu\\
&\Longleftrightarrow & bY_\mu Y_\nu+Z_\mu Z_\nu &= Z_\nu (1+qY_\mu Z^{-1}_\mu Y_\nu Z^{-1}_\nu)Z_\mu\\
&                    &                          &= Z_\mu Z_\nu + q^3 Y_\mu Y_\nu\\
&\Longleftrightarrow &                        b &= q^3.
\end{align*}

In the following we show $b=q^3$ is also sufficient. There are eight cases in Proposition \ref{thm:exchange} to check. For instance, we verify Case 4. By definition, We have
\begin{align*}
F_\tau(X_i) &=Y_\mu Z^{-1}_\mu Y_\nu Z^{-1}_\nu  & F_{\tau'}(X_i) &=Y'_\mu Z'^{-1}_\mu Y'_\nu Z'^{-1}_\nu\\
F_\tau(X_j) &=Z_\mu Z_\nu                        & F_{\tau'}(X_j) &= Y'^{-1}_\mu Y'^{-1}_\nu \\
F_\tau(X_k) &=Y^{-1}_\nu H^s_\zeta               & F_{\tau'}(X_k)&=Z'_\mu H^s_\zeta \\
F_\tau(X_m) &=Y^{-1}_\mu H^t_\eta                & F_{\tau'}(X_m) &=Z'_\nu H^t_\eta
\end{align*}
for some $s,t\in \{0,1,2\}.$

Then we have
\begin{align*}
\widehat{\varphi}_{\mu\nu}\circ F_{\tau'}(X_i)
&=\widehat{\varphi}_{\mu\nu}(Y'_\mu Z'^{-1}_\mu Y'_\nu Z'^{-1}_\nu)\\
&= (bY_\mu Y_\nu+Z_\mu Z_\nu)^{-1}Z_\nu\ \ b^{-1}Y^{-1}_\mu (bY_\mu Y_\nu+Z_\mu Z_\nu)\\
&\ \ \ \ \ \   (bY_\mu Y_\nu+Z_\mu Z_\nu)^{-1}Z_\mu\ \ b^{-1}Y^{-1}_\nu (bY_\mu Y_\nu+Z_\mu Z_\nu) \\
&= b^{-2}  (bY_\mu Y_\nu+Z_\mu Z_\nu)^{-1} Z_\nu Y^{-1}_\mu Z_\mu Y^{-1}_\nu (bY_\mu Y_\nu+Z_\mu Z_\nu)\\
&=  b^{-2}  (bY_\mu Y_\nu+Z_\mu Z_\nu)^{-1} (bq^4Y_\mu Y_\nu+q^4Z_\mu Z_\nu) Z_\nu Y^{-1}_\mu Z_\mu Y^{-1}_\nu \\
&=  b^{-2}q^4  Z_\nu Y^{-1}_\mu Z_\mu Y^{-1}_\nu \\
&=  b^{-2}q^6  Z_\nu Y^{-1}_\nu Z_\mu  Y^{-1}_\mu \\
&=  b^{-2}q^6  F_\tau(X_i^{-1}) \\
&=  b^{-2}q^6  F_\tau\circ \widehat{\Delta}_i(X_i) \\
&=  F_\tau\circ \widehat{\Delta}_i(X_i)
\end{align*}
due to the assumption that $b=q^3.$

And
\begin{align*}
\widehat{\varphi}_{\mu\nu}\circ F_{\tau'}(X_j)
&=\widehat{\varphi}_{\mu\nu}(Y'^{-1}_\mu Y'^{-1}_\nu)\\
&= Z^{-1}_\nu(bY_\mu Y_\nu+Z_\mu Z_\nu)Z^{-1}_\mu(bY_\mu Y_\nu+Z_\mu Z_\nu)\\
&= (b Z^{-1}_\nu Y_\mu Y_\nu Z^{-1}_\mu+1)(bY_\mu Y_\nu+Z_\mu Z_\nu)\\
&= (bq^{-2}  Y_\mu Y_\nu Z^{-1}_\nu Z^{-1}_\mu+1)(b  Y_\mu Y_\nu Z^{-1}_\nu Z^{-1}_\mu+1)Z_\mu Z_\nu \\
&= (q  Y_\mu Y_\nu Z^{-1}_\nu Z^{-1}_\mu+1)(q^3  Y_\mu Y_\nu Z^{-1}_\nu Z^{-1}_\mu+1)Z_\mu Z_\nu \\
&= (1+qF_\tau(X_i))(1+q^3F_\tau(X_i))F_\tau(X_j)\\
&= F_\tau((1+q X_i)(1+q^3 X_i)X_j)\\
&= F_\tau \circ \widehat{\Delta}_i(X_j).
\end{align*}

And
\begin{align*}
\widehat{\varphi}_{\mu\nu}\circ F_{\tau'}(X_k)
&=\widehat{\varphi}_{\mu\nu}(Z'_\mu)H^s\\
&= b(bY_\mu Y_\nu+Z_\mu Z_\nu)^{-1} Y_\mu H^s_\zeta\\
&= b[Y_\mu Y_\nu(b+Y^{-1} _\mu Y^{-1}_\nu Z_\mu Z_\nu)]^{-1}  Y_\mu H^s_\zeta\\
&= b(b+Y^{-1} _\mu Y^{-1}_\nu Z_\mu Z_\nu)^{-1}  Y^{-1}_\nu H^s_\zeta\\
&= b(b+ q^4 Z_\mu Z_\nu Y^{-1} _\mu Y^{-1}_\nu)^{-1} Y^{-1}_\nu H^s_\zeta \\
&= (1+ q Z_\mu Z_\nu Y^{-1} _\mu Y^{-1}_\nu)^{-1} Y^{-1}_\nu H^s_\zeta \\
&= (1+qF_\tau(X_i)^{-1})F_\tau(X_k)\\
&= F_\tau((1+qX_i^{-1})X_k)\\
&= F_\tau \circ \widehat{\Delta}_i(X_k).
\end{align*}

It is same for $X'_m$ due to the symmetry of $\mu,\nu.$
\end{proof}

\begin{theorem}
Suppose the decorated ideal triangulations $\tau$ and $\tau'$ have the underlying ideal triangulations
$\lambda$ and $\lambda'$ respectively.
The following diagram is commutative if and only if $a=q^{-2}, b=q^{3}$.
$$
\begin{CD}
\widehat{\mathcal T}_\lambda^q @>F_\tau>> \widehat{\mathcal K}^q_{\tau}\\
@A\Phi^q_{\lambda,\lambda'} AA @AA\Psi^q_{\tau,\tau'}(a,b) A \\
\widehat{\mathcal T}_{\lambda'}^q @>F_{\tau'}>>  \widehat{\mathcal K}^q_{\tau'}
\end{CD}
$$
\end{theorem}

\begin{proof} By Theorem \ref{thm:Penner1},
$\tau$ and $\tau'$ are connected by a sequence
$\tau=\tau_{(0)}$, $\tau_{(1)}$, \dots,
$\tau_{(n)}=\tau'$ where each $\tau_{(k+1)}$ is obtained
from $\tau_{(k)}$ by a reindexing or a mark rotation or a diagonal exchange. For a reindexing, the diagram is always commutative. By Lemma \ref{thm:a} and \ref{thm:b}, the the diagram is always commutative if and only if $a=q^{-2}, b=q^{3}$.
\end{proof}

Recall that the quantum Teichm\"uller space of $S$ is defined as the algebra
$$
\widehat {\mathcal{T}}^q_S= \bigg(
\bigsqcup_{\lambda\in\Lambda(S)}
\widehat{\mathcal{T}}^q_{\lambda}\bigg)/\sim
$$
where the relation $\sim$ is defined by the property that, for
$X\in \widehat{\mathcal{T}}^q_{\lambda}$ and $X'\in
\widehat{\mathcal{T}}^q_{\lambda'}$,
$$
X \sim X' \Leftrightarrow X=\Phi^q_{\lambda,\lambda'}(X').
$$

And the generalized Kashaev algebra $\widehat{\mathcal{K}}^q_S(a,b)$ associated to a surface $S$ is defined as the algebra
$$
\widehat{\mathcal{K}}^q_S(a,b)= \bigg(
\bigsqcup_{\tau \in \triangle(S)}
\widehat{\mathcal{K}}^q_\tau(a,b)\bigg)/\sim
$$
where the relation $\sim$ is defined by the property that, for
$X\in \widehat{\mathcal{K}}^q_\tau(a,b)$ and $X'\in
\widehat{\mathcal{K}}^q_{\tau'}(a,b)$,
$$
X \sim X' \Leftrightarrow X=\Psi^q_{\tau,\tau'}(a,b)(X').
$$

\begin{corollary}\label{cor:homo} The homomorphism $F_\tau$ induces a homomorphism $\widehat {\mathcal{T}}^q_S \to \widehat{\mathcal{K}}^q_S(a,b)$ if and only if $a=q^{-2}, b=q^3$.
\end{corollary}

\subsection{Quotient algebra}

Furthermore, consider the element
$$H=q^{-\sum_{i<j}\sigma^\lambda_{ij}}X_1X_2...X_{3m}\in \mathcal T_\lambda^q.$$ It is proved in \cite{Liu1}(Proposition 14) that $H$ is independent of the ideal triangulation $\lambda$. Therefore $H$ is a well-defined element of the
quantum Teichm\"uller space $\widehat{\mathcal{T}}^q_S$.

\begin{theorem}\label{thm:homo} The homomorphism $F_\tau$ induces a homomorphism $$\widehat {\mathcal{T}}^q_S/(q^{-2m}H)\to \widehat{\mathcal{K}}^q_S(q^{-2},q^3)$$ where $(q^{-2m}H)$ is the ideal generated by $q^{-2m}H$.
\end{theorem}

\begin{proof} We only need to show that $F_\tau(q^{-2m}H)=1$ for any arbitrary decorated ideal triangulation $\tau$. In fact
$$F_\tau(X_1X_2...X_{3m})=q^{\delta_{\mu_1\nu_1}}H^{s_1}_{\mu_1}H^{t_1}_{\nu_1}...q^{\delta_{\mu_{3m}\nu_{3m}}}H^{s_{3m}}_{\mu_{3m}}H^{t_{3m}}_{\nu_{3m}}.$$
where the edge $\lambda_i$ bounds the $s_i-$side of $\tau_{\mu_i}$ and the $t_i-$side of $\tau_{\nu_i}$ for $i=1,...,3m.$

Since $H^s_\mu$ and $H^t_\nu$ are commutative when $\mu\neq \nu$, we may collect the terms indexed by the same ideal triangle by commutating the terms indexed by different ideal triangles. The right hand side of the above identity is equal to
$$\prod_{\mu=1}^{2m}P_\mu,$$
where $P_\mu$ is the product of terms involving the ideal triangle $\tau_\mu.$

\textbf{Case 1.} If $\tau_\mu$ is embedded, then $P_\mu=H^r_\mu H^s_\mu H^t_\mu,$ where $\{r,s,t\}=\{0,1,2\}$.

When $(r,s,t)$ is an even permutation of $0,1,2$, we have $P_\mu=1.$

Suppose the $r-$side, the $s-$side and the $t-$side of $\tau_\mu$ are bounded by edges $\lambda_i$,  $\lambda_j$ and  $\lambda_k$ respectively. Then $i\leq j \leq k$ since this order is preserved when we commutate the terms indexed by different ideal triangles. Denote by $\sigma^\mu_{ij}$ the number of corners of $\tau_\mu$ delimited by $\lambda_i$ from the left and delimited by $\lambda_j$ from the right minus the number of corners of $\tau_\mu$ delimited by $\lambda_j$ from the left and delimited by $\lambda_i$ from the right. Then
$$\sigma^\mu_{ij}+\sigma^\mu_{jk}+\sigma^\mu_{ik}=-1-1+1=-1.$$
Therefore $$P_\mu=1=q^{1+\sigma^\mu_{ij}+\sigma^\mu_{jk}+\sigma^\mu_{ik}}.$$

When $(r,s,t)$ is an odd permutation of $0,1,2$, we have $P_\mu=q^2.$ And
$$\sigma^\mu_{ij}+\sigma^\mu_{jk}+\sigma^\mu_{ik}=1+1-1=1.$$
Therefore $$P_\mu=q^2=q^{1+\sigma^\mu_{ij}+\sigma^\mu_{jk}+\sigma^\mu_{ik}}.$$

\textbf{Case 2.} If $\tau_\mu$ is not embedded, then $P_\mu=q^{\sigma_{sr}}H^r_\mu H^s_\mu H^t_\mu$ or $P_\mu=q^{\sigma_{ts}}H^r_\mu H^s_\mu H^t_\mu$.
When $(r,s,t)$ is an even permutation of $0,1,2$, we have $P_\mu=q\cdot 1=q.$ When $(r,s,t)$ is an odd permutation of $0,1,2$, we have $P_\mu=q^{-1}\cdot q^2=q.$ And we always have $$\sigma^\mu_{ij}+\sigma^\mu_{jk}+\sigma^\mu_{ik}=0.$$
Therefore
$$P_\mu=q=q^{1+\sigma^\mu_{ij}+\sigma^\mu_{jk}+\sigma^\mu_{ik}}.$$

Combining the two cases, we obtain
$$F_\tau(X_1X_2...X_{3m})=\prod_{\mu=1}^{2m}P_\mu=\prod_{\mu=1}^{2m}q^{1+\sigma^\mu_{ij}+\sigma^\mu_{jk}+\sigma^\mu_{ik}}=q^{2m+\sum_{i<j}\sigma^\lambda_{ij}}.$$
Thus $F_\tau(q^{-2m}H)=1.$
\end{proof}

\section*{Appendix: Kashaev coordinates and Penner coordinates}

We review the relationship of Kasheev coordinate and Penner coordinates following \cite{Kas1} and \cite{Te}.

A \emph{decorated hyperbolic metric} $(d,r)$ on $S,$ introduced by Penner \cite{Pen}, is a complete hyperbolic metric $d$ so that each end is cusp type and each cusp $c_i$ is assigned a positive number $r_i$. The \emph{decorated Teichm\"uller space} is the space of isotopy class of decorated hyperbolic metrics. For each decorated hyperbolic metric $(d,r)$, at each cusp $c_i$, there is a horocycle with boundary length $r_i$. Under a decorated hyperbolic metric, each edge of an ideal triangulation of a punctured surface $S$ is realized as a geodesic running from one puncture to another. Penner coordinate $\delta(e)$ at an edge $e$ is the signed distance between two horocycles bounding cusps $c_i$ and $c_j$ if the edge $e$ runs from $c_i$ to $c_j$. Denote by $\overline{\mathcal{T}}_\lambda$ the decorated Teichm\"uller space parameterized by Penner coordinates associated to the ideal triangulation $\lambda$.

Let $\tau$ be a decorated ideal triangulation with the underlying ideal triangulation $\lambda.$
Let $\mathcal{K}_\tau=\mathbb{R}^{4m}=\{(\ln y_1,\ln z_1,...,\ln y_{2m},\ln z_{2m})\}$ be the space of Kashaev coordinates. There is a map $f:\overline{\mathcal{T}}_\lambda \to \mathcal{K}_\tau$ defined as follows.

For an ideal triangle $\tau_i$ (embedded or not) with a marked corner, there are three sides which correspond to the three half-edges incident to the vertex $\tau_\mu^*$ of the dual graph. The three sides are numerated by $0,1,2$ in the counterclockwise order such that the $0-$side is opposite to the marked corner. Denote by $\lambda_i^0, \lambda_i^1, \lambda_i^2$ the edges (two of them may coincide) bounding the three sides of $\tau_i$. We define
$$y_i=e^{\frac12(\delta(\lambda^i_1)-\delta(\lambda^i_0))}, \ \ \ \ \  z_i=e^{\frac12(\delta(\lambda^i_2)-\delta(\lambda^i_0))}.$$

\begin{prop}[Kashaev \cite{Kas1}]\label{thm:k-exact} The following sequence is exact:
$$1\longrightarrow \mathbb{R}_+ \longrightarrow \overline{\mathcal{T}}_\lambda \xrightarrow{f} \mathcal{K}_\tau
\longrightarrow H^1(S,\mathbb{R}) \longrightarrow 0.$$
\end{prop}

\begin{prop}[Kashaev \cite{Kas1}]\label{thm:k-diff} If $(\ln y_1,\ln z_1,...,\ln y_{2m},\ln z_{2m})=f(\delta(\lambda_1),...,\delta(\lambda_{3m}))$, then the two 2-forms are equal:
$$\sum_{\mu_1}^{2m}d \ln y_\mu \wedge d \ln z_\mu= f^*(\sum_{\mu=1}^{2m}
(d \delta(\lambda_i)\wedge d \delta(\lambda_j) + d \delta(\lambda_j)\wedge d \delta(\lambda_k)
+ d \delta(\lambda_k)\wedge d \delta(\lambda_i)),$$
where $\lambda_i, \lambda_j, \lambda_k$ are edges bounding the three sides of $\tau_\mu$ in the counterclockwise order.
\end{prop}

\begin{prop}[Kashaev \cite{Kas1}]\label{thm:k-com} Suppose that the decorated ideal triangulations $\tau$ and $\tau'$ have the underlying ideal triangulations $\lambda$ and $\lambda'$ respectively. The following diagram is commutative:
$$
\begin{CD}
\overline{\mathcal{T}}_\lambda @>f>> \mathcal{K}_\tau\\
@VVV @VVV \\
\overline{\mathcal{T}}_{\lambda'} @>f>>  \mathcal{K}_{\tau'}
\end{CD}
$$
where the two vertical maps are corresponding coordinate changes. The coordinate changes of Kashaev coordinates are given in Definition \ref{def:coor-change}.

\end{prop}

\begin{figure}[ht!]
\labellist\small\hair 2pt
\pinlabel $\longrightarrow$ at 205 81
\pinlabel $\varphi_{ij}$ at 205 99
\pinlabel $\tau_i$ at 50 85
\pinlabel $\tau_j$ at 107 85
\pinlabel $*$ at 13 85
\pinlabel $*$ at 144 85
\pinlabel $a$ at 31 135
\pinlabel $b$ at 31 33
\pinlabel $c$ at 129 33
\pinlabel $d$ at 129 135
\pinlabel $l$ at 89 111
\pinlabel $\tau'_i$ at 329 107
\pinlabel $\tau'_j$ at 329 48
\pinlabel $*$ at 329 144
\pinlabel $*$ at 329 13
\pinlabel $a$ at 280 135
\pinlabel $b$ at 280 33
\pinlabel $c$ at 376 33
\pinlabel $d$ at 376 135
\pinlabel $m$ at 290 89

\endlabellist
\centering
\includegraphics[scale=0.5]{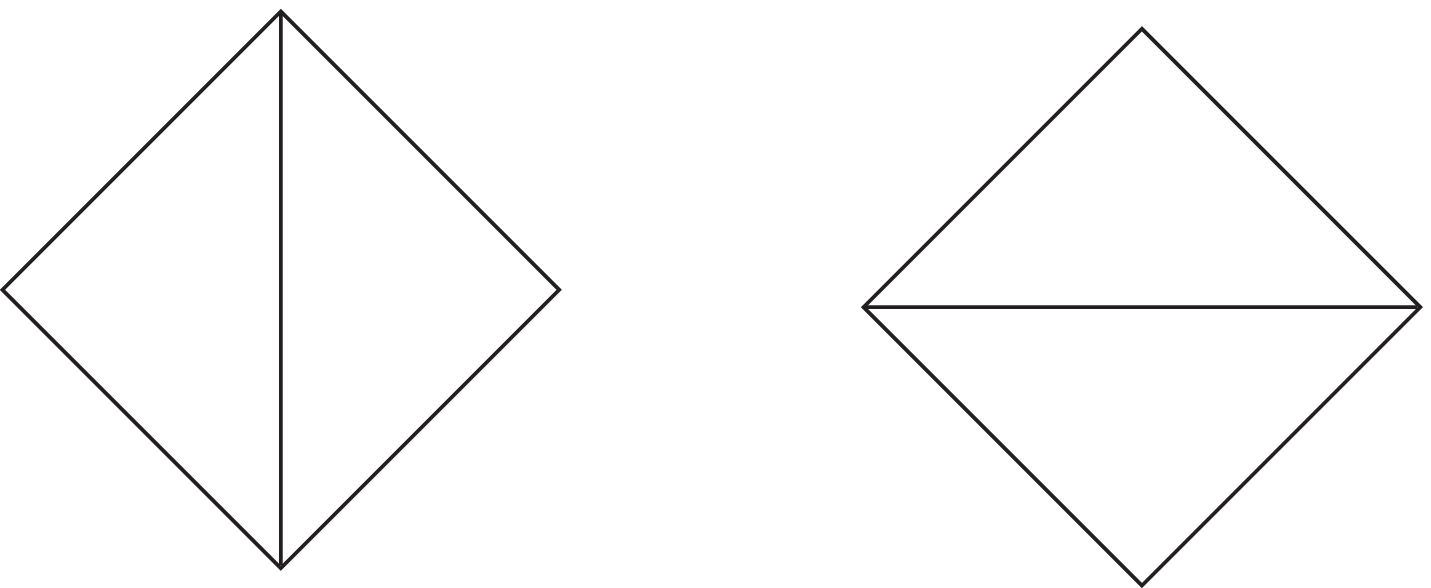}
\caption{}
\label{fig:exchange}
\end{figure}

\begin{proof} For a reindexing, the conclusion is obvious. For a mark rotation, the conclusion is easily proved by applying the definition of $(y_i,z_i).$ For a diagonal exchange, we need to use the famous Ptolemy relation for Penner coordinates.

In Figure \ref{fig:exchange}, denote by $a,b,c,d,l$ and $m$ the Penner coordinates of the corresponding edges. If the ideal triangles are not embedded, some of the numbers $a,b,c,d$ may equal. The Ptolemy relation is
$$e^{\frac12(l+m)}=e^{\frac12(a+c)}+e^{\frac12(b+d)}$$
which holds in spite of whether the ideal triangles $\tau_i,\tau_j$ are embedded or not.

We show the relation between $(y_i,z_i,y_j,z_j)$ and $(y'_i,z'_i,y'_j,z'_j)$ in Definition \ref{def:coor-change} holds. In fact,
\begin{align*}
\frac{z_j}{y_iy_j+z_iz_j}&=\frac{e^{\frac12(d-l)}}{e^{\frac12(a-l)}e^{\frac12(c-l)}+e^{\frac12(b-l)}e^{\frac12(d-l)}}
&(\mbox{by definition})\\
&=\frac{e^{\frac12(d+l)}}{e^{\frac12(a+c)}+e^{\frac12(b+d)}}\\
&=\frac{e^{\frac12(d+l)}}{e^{\frac12(l+m)}} &(\mbox{by Ptolemy relation})\\
&=e^{\frac12(d-m)}\\
&=y'_i.
\end{align*}
The same calculation can be used to verify the formula of $z'_i,y'_j,z'_j.$

\end{proof}

\section*{Acknowledgment}

The authors would like to thank Francis Bonahon and Feng Luo for encouragement and helpful comments, Liang Kong and Hua Bai for helpful discussion. A part of the work of this paper was done when the second author was visiting Chern Institute of Mathematics, Tianjin, China. He would like to take the opportunity to thank Chern Institute for hospitality.

\bibliographystyle{amsplain}

\end{document}